\documentclass[11pt]{article}
\usepackage[a4paper, total={161mm, 241mm}]{geometry}
\usepackage{graphicx}
\usepackage{amssymb}
\usepackage{amsmath}
\usepackage{amsthm}
\usepackage{bm}
\usepackage{tikz}
\usepackage{url}
\usepackage{algorithmic}
\usepackage{algorithm}
\usepackage{soul}
\usepackage{subfigure}
\usepackage{todonotes}
\renewcommand{\ul}[1]{\underline{#1}}
\newcommand{\ddt}{\frac{\text{d}}{\text{d}t}}

\usepackage{hyperref}
\hypersetup{
  pdftitle={A time adaptive multirate Dirichlet-Neumann waveform relaxation method for heterogeneous coupled heat equations},
  pdfauthor={P. Meisrimel, A. Monge, and P. Birken},
  pdfsubject={Mathematics, Numerical Analysis},
  pdfkeywords={Thermal Fluid-Structure Interaction, Coupled Problems, Dirichlet-Neumann Method, Multirate, Time Adaptivity, Waveform Relaxation}
}

\expandafter\def\expandafter\UrlBreaks\expandafter{\UrlBreaks
  \do\a\do\b\do\c\do\d\do\e\do\f\do\g\do\h\do\i\do\j%
  \do\k\do\l\do\m\do\n\do\o\do\p\do\q\do\r\do\s\do\t%
  \do\u\do\v\do\w\do\x\do\y\do\z\do\A\do\B\do\C\do\D%
  \do\E\do\F\do\G\do\H\do\I\do\J\do\K\do\L\do\M\do\N%
  \do\O\do\P\do\Q\do\R\do\S\do\T\do\U\do\V\do\W\do\X%
  \do\Y\do\Z}

\title{A time adaptive multirate Dirichlet--Neumann waveform relaxation method for heterogeneous coupled heat equations}
\author{Peter Meisrimel$^{\mbox{\tiny\rm 1}}$, Azahar Monge$^{\mbox{\tiny\rm 1,2}}$, Philipp Birken$^{\mbox{\tiny\rm 1}}$}
\begin{document}
\maketitle
\baselineskip=0.9
\normalbaselineskip
\vspace{-3pt}
\begin{center}{\footnotesize\em $^{\mbox{\tiny\rm 1}}$Centre for the
    mathematical sciences, Numerical Analysis, Lund University, Lund, Sweden\\ email: peter.meisrimel\symbol{'100}na.lu.se, azahar.sz\symbol{'100}gmail.com, philipp.birken\symbol{'100}na.lu.se\\
    \footnotesize\em $^{\mbox{\tiny\rm 2}}$Chair of Computational Mathematics, University of Deusto, Spain, Bilbao}
\end{center}

\begin{abstract}
We consider partitioned time integration for heterogeneous coupled heat equations. First and second order multirate, as well as time-adaptive Dirichlet-Neumann Waveform relaxation (DNWR) methods are derived. In 1D and for implicit Euler time integration, we analytically determine optimal relaxation parameters for the fully discrete scheme. 

We test the robustness of the relaxation parameters on the second order multirate method in 2D. DNWR is shown to be very robust and consistently yielding fast convergence rates, whereas the closely related Neumann-Neumann Waveform relaxtion (NNWR) method is slower or even diverges. 

The waveform approach naturally allows for different timesteps in the subproblems. In a performance comparison for DNWR, the time-adaptive method dominates the multirate method due to automatically finding suitable stepsize ratios. Overall, we obtain a fast, robust, multirate and time adaptive partitioned solver for unsteady conjugate heat transfer.
\end{abstract}

\medskip\\
{This research is supported by the Swedish e-Science collaboration eSSENCE, which we gratefully acknowledge.}

\section{Introduction}
We consider efficient numerical methods for the partitioned time integration of coupled multiphysics problems. In a partitioned approach different codes for the sub-problems are reused and the coupling is done by a coupling code which calls interface functions of the segregated codes \cite{causin:05, cristiano:11}. These algorithms are currently an active research topic driven by certain multiphysics applications where multiple physical models or multiple simultaneous physical phenomena involve solving coupled systems of partial differential equations (PDEs). 

An example of this is fluid structure interaction (FSI) \cite{vanBrummelen:11, bremicker:17}. Our prime motivation is thermal interaction between fluids and structures, also called conjugate heat transfer. There are two domains with jumps in the material coefficients across the connecting interface. Conjugate heat transfer plays an important role in many applications and its simulation has proved essential \cite{banka:05}. Examples for thermal fluid structure interaction are cooling of gas-turbine blades, thermal anti-icing systems of airplanes \cite{buchli:10}, supersonic reentry of vehicles from space \cite{mehta:05,hinrad:06}, gas quenching, which is an industrial heat treatment of metal workpieces \cite{hefiba:01,stshle:06} or the cooling of rocket nozzles \cite{kohoha:13,kotirh:13}. 

The most common form of coupling is a Dirichlet-Neumann (DN) approach, in which one problem has a Dirichlet boundary condition on the shared interface, while the other one uses a Neumann boundary condition. In the iteration, they provide each other with the suitable boundary data, i.e. a flux or the interface value. Thus, there is a connection to Domain Decomposition methods. 

From the partitioned time integration, we require that it allows for variable and adaptive time steps, preserves the order of the time integration in the subsolvers, and it should be robust and fast. A technique that promises to deliver these properties is the so called Waveform relaxation (WR). An iteration requires solving the subproblems on a time window. Thereby, continuous interface functions, obtained via suitable interpolation, are provided from the respective other problem. WR methods were originally introduced by \cite{lelarasmee:82} for ordinary differential equation (ODE) systems, and used for the first time to solve time dependent PDEs in \cite{gander:98,giladi:02}. They allow the use of different spatial and time discretizations for each subdomain. This is especially useful in problems with strong jumps in the material coefficients \cite{gander:03} or the coupling of different models for the subdomains \cite{Gander:07}. 

A key problem is to make the Waveform iteration fast. A black box approach is to make use of Quasi-Newton methods, leading to Quasi-Newton Waveform iterations \cite{rbmbmb:20}. Here, we follow instead the idea to tailor very fast methods to a specific problem. In particular, we consider the Neumann-Neumann Waveform relaxation (NNWR) and Dirichlet-Neumann Waveform relaxation (DNWR) method of Gander et al. \cite{Fwok:14,gander:16}, which are WR methods based on the classical Neumann-Neumann and Dirichlet-Neumann iterations. The DNWR method is serial, whereas with NNWR, one can solve the subproblems in parallel. Using an optimal relaxation parameter, convergence in two iterations is obtained for the continuous iteration in 1D. In \cite{MongeBirken:multirate}, a fully discrete multirate NNWR method for heterogeneous coupled heat equations is presented. Optimal relaxation parameters are determined for the 1D case. The method was extended to the time adaptive case in \cite{monbirDD25:19}.

However, the NNWR method is extremely sensitive to the choice of the relaxation parameter, leading to a lack of robustness. In this paper, we therefore focus on the DNWR method, see also \cite{Monge2019a}. The standard DN method was known to be a very fast method for thermal interaction between air and steel \cite{biquhm:11,birquihame:10,biglkm:15}. This was thoroughly analyzed for the fully discrete case for two coupled heat equations with different material properties in \cite{Monge:2017} and for coupled Laplace equations in \cite{goebir:20}. Thus, we can expect the waveform variant to be a fast solver not only when using an optimal relaxation parameter. The technique employed here to determine optimal relaxation parameters follows \cite{MongeBirken:multirate}, which considers the fully discrete iteration for a 1D model problem. Then, using the Toeplitz structure of the arising matrices, a formula for the spectral radius of the iteration matrix can be found and the optimal relaxation parameter can be analytically determined. 

We present first and second order multirate WR methods, as well as a second order time adaptive method. The time integration methods we use as a base are implicit Euler and a second order singly diagonally implicit Runge-Kutta (SDIRK2) method. The optimal relaxation parameter $\Theta_{opt}$ from the analysis of 1D and implicit Euler yields good results for 2D and SDIRK2. We show how to adapt $\Theta_{opt}$ for use in the multirate and time-adaptive setting to get good convergence rates. Additionally, we experimentally show that the convergence results also extend to non-square geometries. Overall, our DNWR method yields a fast and robust solver for unsteady conjugate heat transfer.
%
%
\section{Model problem}
%
The unsteady transmission problem reads as follows, where we consider a domain $\Omega \subset \mathbb{R}^d$ which is cut into two subdomains $\Omega = \Omega_1 \cup \Omega_2$ with transmission conditions at the interface $\Gamma = \partial \Omega_1 \cap \partial \Omega_2$:
\begin{align}
\begin{split}\label{EQ PROB MONO}
\alpha_m \frac{\partial u_m(\bm{x},t)}{\partial t} - \nabla \cdot (\lambda_m \nabla u_m(\bm{x},t)) 
&= 0,\quad (\bm{x}, t) \in \Omega_m \times (0, T_f], \, m=1,2,\\
u_m(\bm{x},t) 
&= 0, \quad (\bm{x}, t) \in \partial \Omega_m \backslash \Gamma \times [0, T_f],\\
u_1(\bm{x},t) 
&= u_2(\bm{x},t), \quad (\bm{x}, t) \in \Gamma \times (0, T_f],\\
\lambda_2 \frac{\partial u_2(\bm{x},t)}{\partial \bm{n}_2} 
&= -\lambda_1 \frac{\partial u_1(\bm{x},t)}{\partial \bm{n}_1}, \quad (\bm{x}, t) \in \Gamma \times (0, T_f],\\
u_m(\bm{x},0) 
&= u_m^0(\bm{x}), \quad \bm{x} \in \Omega_m,
\end{split}
\end{align}
where $\bm{n}_m$ is the outward normal to $\Omega_m$ for $m=1,2$.

The constants $\lambda_1$ and $\lambda_2$ describe the thermal conductivities of the materials on $\Omega_1$ and $\Omega_2$ respectively. $D_1$ and $D_2$ represent the thermal diffusivities of the materials and are defined by
\begin{align*}
D_m = \lambda_m /\alpha_m, \quad \mbox{with} \quad \alpha_m = \rho_m c_{p_m},
\end{align*}
where $\rho_m$ is the density and $c_{p_m}$ the specific heat capacity of the material placed in $\Omega_m$, $m=1,2$.
%
\section{The Dirichlet-Neumann Waveform Relaxation algorithm}
%
The Dirichlet-Neumann Waveform relaxation (DNWR) method is inspired by substructuring methods from Domain Decomposition. The PDEs are solved sequentially using a Dirichlet- respectively Neumann boundary condition with data given from the solution of the other problem, c.f. \cite{Mandal:14,mandalphd}. 

Given an interface solution $g^{(k)} (\bm{x},t)$, $(\bm{x}, t) \in \Gamma \times [0, T_f]$, it consists of the following three-step iteration. Imposing continuity of the solution across the interface, one first finds the local solution $u_1^{(k+1)}(\bm{x},t)$ on $(\bm{x},t) \in \Omega_1 \times [0, T_f]$ by solving the Dirichlet problem:
\begin{align}
\begin{split}\label{EQ CONT DIR PROB}
\alpha_1 \frac{\partial u_1^{(k+1)}(\bm{x},t)}{\partial t} - \nabla \cdot (\lambda_1 \nabla u_1^{(k+1)}(\bm{x},t)) = 0, \quad (\bm{x}, t) \in \Omega_1 \times (0, T_f],\\
u_1^{(k+1)}(\bm{x},t) = 0, \quad (\bm{x}, t) \in \partial \Omega_1 \backslash \Gamma \times [0, T_f], \\
u_1^{(k+1)}(\bm{x},t) = g^{(k)}(\bm{x},t), \quad (\bm{x}, t) \in \Gamma \times [0, T_f], \\
u_1^{(k+1)}(\bm{x},0) = u_1^{0}(\bm{x}), \quad \bm{x} \in \Omega_1.
\end{split}
\end{align}
The typical initial guess is $g^{(0)}(\bm{x}, t) = u(\bm{x}, 0)\big|_\Gamma$, i.e., extrapolation of the interface initial value.

Then, imposing continuity of the heat fluxes across the interface, one finds the local solution $u_2^{(k+1)}(\bm{x},t)$ on $\Omega_2$ by solving the Neumann problem:
\begin{align}
\begin{split}\label{EQ CONT NEU PROB}
\alpha_2 \frac{\partial u_2^{(k+1)}(\bm{x},t)}{\partial t} - \nabla \cdot (\lambda_2 \nabla u_2^{(k+1)}(\bm{x},t)) = 0, \quad (\bm{x}, t) \in \Omega_2 \times (0, T_f],\\
u_2^{(k+1)}(\bm{x},t) = 0, \quad (\bm{x}, t) \in \partial \Omega_2 \backslash \Gamma \times [0, T_f], \\
\lambda_2 \frac{\partial u_2^{(k+1)}(\bm{x},t)}{\partial \bm{n}_2} = -\lambda_1 \frac{\partial u_1^{(k+1)}(\bm{x},t)}{\partial \bm{n}_1}, \quad (\bm{x},t) \in \Gamma \times [0, T_f], \\
u_2^{(k+1)}(\bm{x},0) = u_2^0(\bm{x}), \quad \bm{x} \in \Omega_2.
\end{split}
\end{align}
Finally, the interface values are updated with
\begin{align}
\label{EQ CONT UPDATE PROB}
g^{(k+1)}(\bm{x}, t) = \Theta u_2^{(k+1)}(\bm{x}, t) + (1 - \Theta) g^{(k)}(\bm{x}, t), \quad (\bm{x},t) \in \Gamma \times [0, T_f],
\end{align}
where $\Theta \in (0,1]$ is the relaxation parameter. Note that choosing an appropriate relaxation parameter is crucial to get a fast convergence rate. In \cite{gander:16}, the optimal relaxation parameter in 1D has been proven to be $\Theta = 1/2$ for $\lambda_1=\lambda_2=\alpha_1=\alpha_2=1$ and subdomains of equal size. If one uses the optimal relaxation parameter, two iterations are enough for subdomains of equal size.
%
\section{Semidiscrete method}
\begin{figure}
\begin{center}
\begin{tikzpicture}[scale = 1.25]
\draw [line width=0.5mm] (0,0) -- (8,0) -- (8,4) -- (0,4) -- (0,0); 
\draw [line width=0.5mm] (4,0) -- (4,4); 

\draw (1, 0) -- (1, 4); \draw (2, 0) -- (2, 4); \draw (3, 0) -- (3, 4);
\draw (5, 0) -- (5, 4); \draw (6, 0) -- (6, 4); \draw (7, 0) -- (7, 4);
\draw (0, 1) -- (8, 1); \draw (0, 2) -- (8, 2); \draw (0, 3) -- (8, 3);
\draw (0, 1) -- (1, 0); \draw (0, 2) -- (2, 0); \draw (0, 3) -- (3, 0);
\draw (0, 4) -- (4, 0); \draw (1, 4) -- (5, 0); \draw (2, 4) -- (6, 0);
\draw (3, 4) -- (7, 0); \draw (4, 4) -- (8, 0); \draw (5, 4) -- (8, 1);
\draw (6, 4) -- (8, 2); \draw (7, 4) -- (8, 3);

\node at (2, -0.25) {\color{red}\large$\Omega_1$};
\node at (4, -0.25) {\color{green}\large$\Gamma$};
\node at (6, -0.25) {\color{blue}\large$\Omega_2$};
%
\foreach \x in {1,2,3}
    \foreach \y in {1,2,3}
        \draw [color = red, line width = 0.6mm] (\x,  \y) circle(3pt);

\foreach \y in {1,2,3}
    \draw [color = green, line width = 0.6mm] (4,  \y) circle(3pt);

\foreach \x in {5,6,7}
    \foreach \y in {1,2,3}
        \draw [color = blue, line width = 0.6mm] (\x,  \y) circle(3pt);
\end{tikzpicture}
\end{center}
\caption{Splitting of $\Omega$ and finite element triangulation.}
\label{FIG DOMAIN FE}
\end{figure}
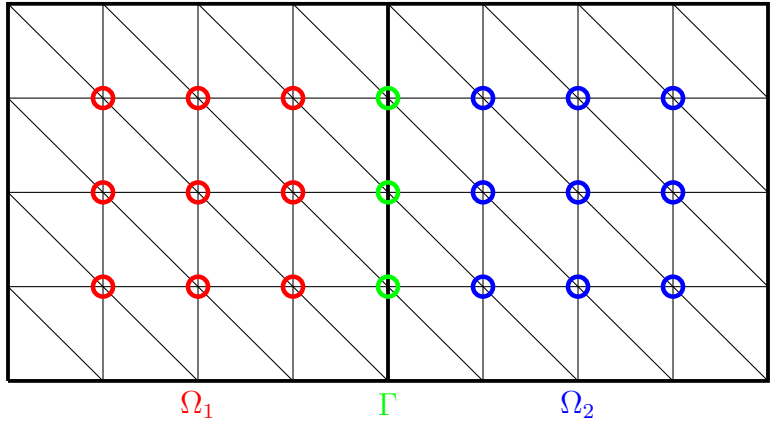

%
We now describe a rather general space discretization of \eqref{EQ CONT DIR PROB}-\eqref{EQ CONT UPDATE PROB}. The core property we need is that the meshes of $\Omega_1$ and $\Omega_2$ share the same nodes on $\Gamma$ as shown in Figure \ref{FIG DOMAIN FE}. Furthermore, we assume that there is a specific set of unknowns associated with the interface nodes. Otherwise, we allow at this point for arbitrary meshes on both sides. 
Then, letting $\bm{u}_I^{(m)}: [0, T_f] \rightarrow \mathbb{R}^{R_m}$ where $R_m$ is the number of grid points on $\Omega_m$, $m=1,2$, and $\bm{u}_{\Gamma}: [0, T_f] \rightarrow \mathbb{R}^{s}$, where $s$ is the number of grid points at the interface $\Gamma$, we can write a general discretization of the first two equations in \eqref{EQ CONT DIR PROB} and \eqref{EQ CONT NEU PROB}, respectively, in a compact form as:
\begin{align}
\label{EQ SEMI DISCR DIR}
\bm{M}_{II}^{(1)} \dot{\bm{u}}_I^{(1),(k+1)}(t) 
+ \bm{A}_{II}^{(1)} \bm{u}_I^{(1),(k+1)}(t) 
= 
-\bm{M}_{I \Gamma}^{(1)} \dot{\bm{u}}_{\Gamma}^{(k)}(t) 
- \bm{A}_{I \Gamma}^{(1)} \bm{u}_{\Gamma}^{(k)}(t), \quad t \in [0, T_f],\\
\label{EQ SEMI DISCR NEU SINGLE}
\bm{M}_{II}^{(2)} \dot{\bm{u}}_I^{(2),(k+1)}(t) 
+ \bm{A}_{II}^{(2)} \bm{u}_I^{(2),(k+1)}(t) 
+ \bm{M}_{I \Gamma}^{(2)} \dot{\bm{u}}_{\Gamma}^{(k+1)}(t) 
+ \bm{A}_{I \Gamma}^{(2)} \bm{u}_{\Gamma}^{(k+1)}(t) 
= \mathbf{0}, \quad t \in [0, T_f],
\end{align}
with initial conditions $\bm{u}_I^{(m)}(0) \in \mathbb{R}^{S_m}$ resp. $\bm{u}_{\Gamma}(0) \in \mathbb{R}^{s}$, $m=1,2$.

To close the system, we need an approximation of the normal derivatives on $\Gamma$. Letting $\phi_j$ be a nodal FE basis function on $\Omega_m$ for a node on $\Gamma$ we observe that the normal derivative of $u_m$ with respect to the interface can be written as a linear functional using Green's formula \cite[pp. 3]{Toselli:2004qr}. Thus, the approximation of the normal derivative is given by
\begin{align*}
\begin{split}
& \lambda_m \int_{\Gamma} \frac{\partial u_m(\bm{x},t)}{\partial \bm{n}_m} \phi_j(\bm{x}) dS 
= \lambda_m \int_{\Omega_m} (\Delta u_m(\bm{x},t) \phi_j(\bm{x}) + \nabla u_m(\bm{x},t) \nabla \phi_j(\bm{x})) d \bm{x} \\
& = \alpha_m \int_{\Omega_m} \frac{d}{dt} u_m(\bm{x},t) \phi_j(\bm{x}) d\bm{x} + \lambda_m \int_{\Omega_m} \nabla u_m(\bm{x},t) \nabla \phi_j(\bm{x}) d \bm{x},\quad t \in [0, T_f], \quad m=1,2.
\end{split}
\end{align*}
Consequently, the equation
\begin{align}
\label{EQ SEMI DISCR NEU GAMMA}
\begin{split}
& 
\bm{M}_{\Gamma I}^{(2)} \dot{\bm{u}}_{I}^{(2),(k+1)}(t) 
+ \bm{A}_{\Gamma I}^{(2)} \bm{u}_{I}^{(2),(k+1)}(t)
+ \bm{M}_{\Gamma \Gamma}^{(2)} \dot{\bm{u}}_{\Gamma}^{(k+1)}(t) 
+ \bm{A}_{\Gamma \Gamma}^{(2)} \bm{u}_{\Gamma}^{(k+1)}(t) 
\\
=&\,\ - \left( 
\bm{M}_{\Gamma I}^{(1)} \dot{\bm{u}}_{I}^{(1),(k+1)}(t)
+ \bm{A}_{\Gamma I}^{(1)} \bm{u}_{I}^{(1),(k+1)}(t) 
+ \bm{M}_{\Gamma \Gamma}^{(1)} \dot{\bm{u}}_{\Gamma}^{(k)}(t) 
+ \bm{A}_{\Gamma \Gamma}^{(1)} \bm{u}_{\Gamma}^{(k)}(t) 
\right), \quad t \in [0, T_f],
\end{split}
\end{align}
is a semi-discrete version of the third equation in \eqref{EQ CONT NEU PROB} and completes the system \eqref{EQ SEMI DISCR DIR}-\eqref{EQ SEMI DISCR NEU SINGLE}. 

Omitting the iteration indices, the system of IVPs defined by \eqref{EQ SEMI DISCR DIR}, \eqref{EQ SEMI DISCR NEU SINGLE} and \eqref{EQ SEMI DISCR NEU GAMMA} is a semidiscretization of \eqref{EQ PROB MONO}. We refer to it as the (semidiscrete) monolithic system and its solution as the monolithic solution.

Then, the semidiscrete DNWR algorithm is as follows: In each iteration $k$, one first solves the Dirichlet problem \eqref{EQ SEMI DISCR DIR}, obtaining $\bm{u}_I^{(1),(k+1)}$. Then, using the function of unknowns $\bm{u}^{(k+1)} = \left( {\bm{u}_I^{(2),(k+1)}}^T {\bm{u}_{\Gamma}^{(k+1)}}^T \right)^T$, one solves the following Neumann problem that corresponds to equations \eqref{EQ SEMI DISCR NEU SINGLE} and \eqref{EQ SEMI DISCR NEU GAMMA}:
\begin{align}
\label{EQ SEMI DISCR NEU}
\bm{M} \dot{\bm{u}}^{(k+1)}(t) + \bm{A} \bm{u}^{(k+1)}(t) = \bm{b}^{(k)}(t), \quad t \in [0, T_f], \quad \bm{u}(0) = \bm{u}_0,
\end{align}
where
\begin{align*}
\bm{M} = \left( \begin{array}{cc}
\bm{M}_{II}^{(2)} & \bm{M}_{I \Gamma}^{(2)} \\
\bm{M}_{\Gamma I}^{(2)} & \bm{M}_{\Gamma \Gamma}^{(2)}\\
\end{array} \right), \quad \bm{A} = \left( \begin{array}{cc}
\bm{A}_{II}^{(2)} & \bm{A}_{I \Gamma}^{(2)} \\
\bm{A}_{\Gamma I}^{(2)} & \bm{A}_{\Gamma \Gamma}^{(2)}\\
\end{array} \right), \quad \bm{b}^{(k)}(t) = \left( \begin{array}{c}
\bm{0} \\
-\bm{q}^{(k+1)}(t)\\
\end{array} \right),
\end{align*}
with the heat flux
\begin{equation}\label{EQ SEMI DISCR HEAT FLUX}
\bm{q}^{(k+1)}(t) = 
\bm{M}_{\Gamma I}^{(1)} \dot{\bm{u}}_{I}^{(1),(k+1)}(t)
+ \bm{A}_{\Gamma I}^{(1)} \bm{u}_{I}^{(1),(k+1)}(t)
+ \bm{M}_{\Gamma \Gamma}^{(1)} \dot{\bm{u}}_{\Gamma}^{(k)}(t)
+ \bm{A}_{\Gamma \Gamma}^{(1)} \bm{u}_{\Gamma}^{(k)}(t)
.
\end{equation}
Finally, the interface solutions are updated by
\begin{align}
\label{EQ SEMI DISCR UPDATE}
\bm{u}_{\Gamma}^{(k+1)}
\gets
\Theta \bm{u}_{\Gamma}^{(k+1)} + (1 - \Theta) \bm{u}_{\Gamma}^{(k)}.
\end{align}
The iteration starts with a function initial guess $\bm{u}_{\Gamma}^{(0)}$, we use $\bm{u}_{\Gamma}^{(0)} \equiv \bm{u}_{\Gamma}(0)$. Since the iteration is done on functions, one would like to terminate when $ \| \bm{u}_{\Gamma}^{(k+1)} - \bm{u}_{\Gamma}^{(k)} \| \leq TOL_{WR} $ is met, where $TOL_{WR}$ is a user defined tolerance. Since we expect the time integration error to grow with $t$, we only compare the update at $T_f$. Our termination criterion is
\begin{equation}
\label{EQ TERMINATION CRIT}
\| \bm{u}_{\Gamma}^{(k+1)}(T_f) - \bm{u}_{\Gamma}^{(k)}(T_f) \|_{\Gamma} \leq TOL_{WR} \cdot \| \bm{u}_\Gamma(0)\|_{\Gamma},
\end{equation}
i.e., the relative update w.r.t. the initial value at the interface. We use the discrete $\mathcal{L}^2$ interface norm, given by
\begin{equation*}
\| \cdot \|_{\Gamma} = \| \cdot \|_2 \Delta x^{(d-1)/2},
\end{equation*}
for a mesh uniform at the interface. Here, $d$ is the spatial dimension of \eqref{EQ PROB MONO}.
%
\section{Multirate time discretization}\label{SEC TIME DISCR}
%
In the case of non-matching material parameters $\alpha_m$, $\lambda_m$, the Dirichlet and the Neumann problems \eqref{EQ SEMI DISCR DIR} and \eqref{EQ SEMI DISCR NEU} have different needs for time-discretization. Consequently, we want the possibility to use different time integration methods and different step-sizes on each subdomain. To this end, we first define interpolation functions in time at the interface. We then present two fully discrete DNWR methods, which use the implicit Euler, resp. the second order singly diagonally implicit Runge-Kutta method (SDIRK2) for time integration. We formulate these for general step-sizes $\Delta t_n$. One obtains the multirate resp. adaptive algorithm by choosing different step-sizes on the different subdomains.

\subsection{Interpolation}
We denote the fully discrete interface solutions resp. heat fluxes by 
\begin{equation*}
\ul{\bm{u}}_{\Gamma}^{(k)} := \{ \bm{u}_\Gamma^{(k), n}\}_{n = 0, \ldots, N_2^{(k)}}
\quad \text{and} \quad
\ul{\bm{q}}^{(k)} := \{ \bm{q}^{(k),n}\}_{n = 0, \ldots, N_1^{(k)}},
\end{equation*}
where $N_m^{(k)}$ is the number of timesteps on $\Omega_m$ in the $k$-th iteration. We define the interpolants as follows:
\begin{equation*}
\mathcal{I}(\ul{\bm{u}}_\Gamma^{(k)}) \in \mathcal{C}\left([0, T_f]; \mathbb{R}^{d_s}\right)
\quad \text{and} \quad
\mathcal{I}(\ul{\bm{q}}^{(k)}) \in \mathcal{C}\left([0, T_f]; \mathbb{R}^{d_s}\right),
\end{equation*}
omitting the dependencies of the interpolants on the time-grids for readability. Now, when deriving the fully discrete DNWR methods, we replace all evaluations of the interface temperatures in the Dirichlet problem \eqref{EQ SEMI DISCR DIR} by evaluations of the interpolant $\mathcal{I}(\ul{\bm{u}}_\Gamma^{(k)})$. 
Analogously, we replace heat flux evaluations for the Neumann problem \eqref{EQ SEMI DISCR NEU} by evaluations of $\mathcal{I}(\ul{\bm{q}}^{(k+1)})$. 

The interpolation is done by a coupling code, such that a solver for a subdomain need not know about the other time grid. Here, we consider \textit{linear polynomial interpolation}.

For readability, we omit dependencies of the time-points on $(k)$ in the following derivations. In case of time-grids varying with $k$, all time-point evaluations correspond to the time-grid of the current iteration. The time-grid of the previous iteration is only used to define the underlying interpolation data.
%
\subsection{Implicit Euler}
%
The implicit Euler method is defined by approximating the derivative term via a standard backward differences. We will use this approximation for derivative terms, i.e.,
\begin{align*}
\ddt\mathcal{I}(\ul{\bm{u}}_\Gamma^{(k)})(t_{n+1}) 
&\approx
\left( \mathcal{I}(\ul{\bm{u}}_\Gamma^{(k)})(t_{n+1}) - \mathcal{I}(\ul{\bm{u}}_\Gamma^{(k)})(t_{n})\right)/\Delta t_n
\\
\text{and} \quad
\dot{\bm{u}}_I^{(1),(k+1)}(t_{n+1})
&\approx
\left( \bm{u}_I^{(1),(k+1)}(t_{n+1}) - \bm{u}_I^{(1),(k+1)}(t_{n})\right)/\Delta t_n.
\end{align*}
In each iteration $k$ of the WR algorithm, given the initial value $\bm{u}_I^{(1),(k+1),0} = \bm{u}_I^{(1)}(0)$ and the function $\mathcal{I}(\ul{\bm{u}}^{(k)}_\Gamma)(t)$, we apply the implicit Euler scheme to the Dirichlet problem \eqref{EQ SEMI DISCR DIR}, which yields
\begin{align}
\label{EQ IE DIR}
\begin{split}
&\left( \bm{M}_{II}^{(1)} + \Delta t\bm{A}_{II}^{(1)} \right) \bm{u}_I^{(1),(k+1),n+1} 
\\ 
=& \,\,
\bm{M}_{II}^{(1)} \bm{u}_I^{(1),(k+1),n} 
- \left( \bm{M}_{I \Gamma}^{(1)} + \Delta t_n \bm{A}_{I \Gamma}^{(1)} \right)\mathcal{I}(\ul{\bm{u}}_\Gamma^{(k)})(t_{n+1})
+ \bm{M}_{I \Gamma}^{(1)} \mathcal{I}(\ul{\bm{u}}_\Gamma^{(k)})(t_{n}), \quad n = 0,\ldots, N_1^{(k+1)} - 1.
\end{split}
\end{align}
We compute the discrete heat fluxes as output of the Dirichlet solver by approximating the derivatives in \eqref{EQ SEMI DISCR HEAT FLUX} using standard backward differences, this yields
\begin{align}\label{EQ IE FLUX}
\begin{split}
\bm{q}^{(k+1), n+1} 
=&\,\,
\left(\bm{M}_{\Gamma I}^{(1)}/\Delta t_n + \bm{A}_{\Gamma I}^{(1)} \right)\bm{u}_I^{(1),(k+1), n+1} 
- \bm{M}_{\Gamma I}^{(1)}/\Delta t_n \bm{u}_I^{(1),(k+1),n} \\
& + \left(\bm{M}_{\Gamma \Gamma}^{(1)}/\Delta t_n + \bm{A}_{\Gamma \Gamma}^{(1)} \right) \mathcal{I}(\ul{\bm{u}}_\Gamma^{(k)})(t_{n+1}) 
- \bm{M}_{\Gamma \Gamma}^{(1)}/\Delta t_n \mathcal{I}(\ul{\bm{u}}_\Gamma^{(k)})(t_{n})
, \quad n = 0,\ldots, N_1^{(k+1)} - 1.
\end{split}
\end{align}
The initial flux $\bm{q}^{(k+1), 0}$, which is required in the interpolation, is analogously computed using standard forward differences to approximate derivative terms.

Next, we rewrite the Neumann problem \eqref{EQ SEMI DISCR NEU} in terms of the vector of unknowns
\begin{equation*}
\bm{u}^{(k+1),n+1} := \Big( {\bm{u}_I^{(2),(k+1),n+1}}^T {\bm{u}_{\Gamma}^{(k+1),n+1}}^T \Big)^T.
\end{equation*}
With $\bm{u}^{(k+1),0} = \bm{u}(0)$ and $\mathcal{I}(\ul{\bm{q}}^{(k+1)})$, the implicit Euler scheme for \eqref{EQ SEMI DISCR NEU} is
\begin{align}\label{EQ IE NEU}
\left( \bm{M} + \Delta t\bm{A} \right) \bm{u}^{(k+1),n+1} 
= \bm{M} \bm{u}^{(k+1),n}
- \Delta t \begin{pmatrix} \bm{0} \\ \mathcal{I}(\ul{\bm{q}}^{(k+1)})(t_{n+1}) \end{pmatrix}
, \quad n = 0,\ldots, N_2^{(k+1)} - 1.
\end{align}
Finally, the interfaces values are updated by
\begin{align*}
\bm{u}_{\Gamma}^{(k+1), n} \gets \Theta \bm{u}_{\Gamma}^{(k+1), n} + (1 - \Theta) \mathcal{I}(\ul{\bm{u}}_{\Gamma}^{(k)})(t_n), \quad n = 0, \ldots, N_2^{(k+1)}.
\end{align*}
Algorithm \ref{ALG IE FULL} shows a pseudocode of this method.
\begin{algorithm}[ht!]
\caption{Pseudocode of the DNWR IE method. On domain $\Omega_m$ we do $N_m$ timesteps of size $\Delta t_m = T_f/N_m$.}
\label{ALG IE FULL}
\begin{algorithmic}
\STATE{\textbf{DNWR\_IE}($T_f$, $N_1$, $N_2$, $(\bm{u}_0^{(1)}, \bm{u}_0^{(2)}, \bm{u}_\Gamma(0))$, $\Theta$, $TOL_{WR}$, $k_{\max}$):}
\STATE{$\ul{\bm{u}}_\Gamma^{(0)} \equiv \bm{u}_{\Gamma} (0)$ Initial guess}
\FOR{$k = 0, \ldots, k_{\max} - 1$}
    \STATE{$\mathcal{I}(\ul{\bm{u}}^{(k)}_\Gamma) \gets \text{Interpolation}(\ul{\bm{u}}_\Gamma^{(k)})$}
    \STATE{$\ul{\bm{q}}^{(k+1)}$ $\gets$ \texttt{SolveDirichlet}($T_f$, $N_1$, $\bm{u}_0^{(1)}$, $\mathcal{I}(\ul{\bm{u}}^{(k)})$)}
    \STATE{$\mathcal{I}(\ul{\bm{q}}^{(k+1)}) \gets \text{Interpolation}(\ul{\bm{q}}^{(q)})$}
    \STATE{$\ul{\bm{u}}_{\Gamma}^{(k+1)}$ $\gets$ \texttt{SolveNeumann}($T_f$, $N_2$, $(\bm{u}_0^{(2)}, \bm{u}_0(x_\Gamma))$, $\mathcal{I}(\ul{\bm{q}}^{(k+1)})$)}
    \STATE{$\ul{\bm{u}}_{\Gamma}^{(k+1)} \gets \Theta \ul{\bm{u}}_{\Gamma}^{(k+1)} + (1 - \Theta) \ul{\bm{u}}_{\Gamma}^{(k)}$}
    \IF{$\| \bm{u}_{\Gamma}^{(k+1), N_2} - \bm{u}_{\Gamma}^{(k), N_2}\|_{\Gamma} < TOL_{WR} \, \| \bm{u}_\Gamma(0)\|_{\Gamma}$}
        \STATE{\textbf{break}}
    \ENDIF
\ENDFOR
\end{algorithmic}
\end{algorithm}
%
\subsection{SDIRK2}\label{SEC SDIRK2}
%
We now introduce a higher order version of the same multirate algorithm. Specifically, we consider the second order singly diagonally implicit Runge-Kutta (SDIRK2) method as a basis to discretize the systems \eqref{EQ SEMI DISCR DIR}, \eqref{EQ SEMI DISCR NEU} and \eqref{EQ SEMI DISCR UPDATE} in time. For a general IVP $\dot{\bm{u}}(t) = \bm{f}(t, \bm{u}(t))$, $\bm{u}(0)= \bm{u}_0$, $t \in [0, T_f]$, SDIRK2 is defined as follows:
\begin{align*}
\begin{split}
\bm{U}_1 
&= \bm{u}_n + a \Delta t_n \bm{f}(t_n + a \Delta t_n, \bm{U}_1), \\
\bm{u}_{n+1} = \bm{U}_2 
&= \bm{u}_n + (1 - a) \Delta t_n \bm{f}(t_n + a \Delta t_n, \bm{U}_1) + a \Delta t_n \bm{f}(t_n + \Delta t_n, \bm{U}_2),
\end{split}
\end{align*}
with $a = 1 - \frac{1}{2} \sqrt{2}$. Here, the so-called stage derivatives are
\begin{align*}
\bm{k}_1 = \bm{f}(t_n + a \Delta t_n, \bm{U}_1)
\approx \dot{\bm{u}}(t_n + a \Delta t_n)
\quad \text{and} \quad
\bm{k}_2 = \bm{f}(t_n + \Delta t_n, \bm{U}_2)
\approx \dot{\bm{u}}(t_n + \Delta t_n).
\end{align*}
In the following we use $\bm{U}_j^{(m)}, \bm{k}_j^{(m)}$, $j = 1, 2$, to denote the stage solutions resp. derivatives on $\Omega_m$, $m = 1,2$.

Using the SDIRK2 scheme to solve the Dirichlet problem \eqref{EQ SEMI DISCR DIR}, with initial value $\bm{u}_I^{(1),(k+1),0} = \bm{u}_I^{(1)}(0)$ and approximating $\bm{U}_1^{(2)}\big|_\Gamma$ by $\mathcal{I}(\ul{\bm{u}}_\Gamma^{(k)})(t_n + a \Delta t_n)$, first yields
\begin{align}\label{EQ SDIRK2 DIR 1}
\begin{split}
& \left( \bm{M}_{II}^{(1)} + a \Delta t \bm{A}_{II}^{(1)} \right) \bm{U}_1^{(1)} \\
= & \,\, \bm{M}_{II}^{(1)} \bm{u}^{(1),(k+1), n} 
- a \Delta t \left(\bm{M}_{I \Gamma}^{(1)} \ddt\mathcal{I}(\ul{\bm{u}}_{\Gamma}^{(k)})(t_n + a \Delta t_n) 
+ \bm{A}_{I \Gamma}^{(1)} \mathcal{I}(\ul{\bm{u}}_{\Gamma}^{(k)})(t_n + a \Delta t)\right).
\end{split}
\end{align}
We then compute the stage heat flux
\begin{align}\label{EQ SDIRK2 FLUX 1}
\bm{q}^{(k+1), n+1}_1
=
\bm{M}_{\Gamma I}^{(1)} \bm{k}_1^{(1)} 
+ \bm{A}_{\Gamma I}^{(1)} \bm{U}_1^{(1)}
+ \bm{M}_{\Gamma \Gamma}^{(1)} \ddt \mathcal{I}(\ul{\bm{u}}_\Gamma^{(k)})(t_{n} + a \Delta t_n)
+ \bm{A}_{\Gamma I}^{(1)} \mathcal{I}(\ul{\bm{u}}_\Gamma^{(k)})(t_{n} + a \Delta t_n).
\end{align}
In both \eqref{EQ SDIRK2 DIR 1} and \eqref{EQ SDIRK2 FLUX 1}, we approximate the derivative terms by
\begin{align*}
\ddt \mathcal{I}(\ul{\bm{u}}_\Gamma^{(k)}) (t_n + a \Delta t_n) 
\approx 
\left( \mathcal{I}(\ul{\bm{u}}_\Gamma^{(k)}) (t_n + a \Delta t_n) - \mathcal{I}(\ul{\bm{u}}_\Gamma^{(k)}) (t_n)\right)/(a \Delta t_n).
\end{align*}
For the second SDIRK2 stage, one solves
\begin{align*}
\begin{split}
&\left( \bm{M}_{II}^{(1)} + a \Delta t \bm{A}_{II}^{(1)} \right) \bm{u}_I^{(1), (k+1), n+1} \\
&= \bm{M}_{II}^{(1)} \left( \bm{u}^{(1),(k+1), n} + (1-a) \Delta t_n\bm{k}_1^{(1)} \right)\\
& \quad - a \Delta t \left(\bm{M}_{I \Gamma}^{(1)} \ddt\mathcal{I}(\ul{\bm{u}}_{\Gamma}^{(k)}) (t_n + \Delta t_n) 
+ \bm{A}_{I \Gamma}^{(1)} \mathcal{I}(\ul{\bm{u}}_{\Gamma}^{(k)})(t_n + \Delta t)\right),
\end{split}
\end{align*}
and computes the second heat flux
\begin{align*}
\begin{split}
\bm{q}^{(k+1), n+1}_2
=&\,\,
\bm{M}_{\Gamma I}^{(1)} \bm{k}_2^{(1)} 
+ \bm{A}_{\Gamma I}^{(1)} \bm{u}_I^{(1), (k+1), n+1} \\ 
& + \bm{M}_{\Gamma \Gamma}^{(1)} \ddt \mathcal{I}(\ul{\bm{u}}_\Gamma^{(k)})(t_{n} + \Delta t_n)
+ \bm{A}_{\Gamma I}^{(1)} \mathcal{I}(\ul{\bm{u}}_\Gamma^{(k)})(t_{n} + \Delta t_n).
\end{split}
\end{align*}
In the computation of the second stage, we use the derivative approximation
\begin{align*}
\begin{split}
& \ddt \mathcal{I}(\ul{\bm{u}}_\Gamma^{(k)}) (t_n + \Delta t_n) 
\\
&\approx 
\left( \mathcal{I}(\ul{\bm{u}}_\Gamma^{(k)}) (t_n + \Delta t_n) 
- \left( \mathcal{I}(\ul{\bm{u}}_\Gamma^{(k)})(t_n) + (1-a) \Delta t_n \ddt \mathcal{I}(\ul{\bm{u}}_\Gamma^{(k)}) (t_n + a \Delta t_n) \right)\right)/(a \Delta t_n).
\end{split}
\end{align*}
This procedure is repeated for all $n = 0, \ldots, N_1^{(k+1)}-1$. The time-point evaluations in our derivative approximations coincide with the time-points for which the SDIRK2 scheme provides discrete solutions, in case of matching time-grids on both subdomains.

We now construct the two interpolants $\mathcal{I}(\ul{\bm{q}}_1^{(k+1)})$ and $\mathcal{I}(\ul{\bm{q}}_2^{(k+1)})$, where the former corresponds to the time-points $0, a \Delta t_0, t_1 + a \Delta t_1, \ldots$. We compute the initial flux $\bm{q}^{(k+1), 0}$, which we use in both flux interpolants, using the second order forward differences:
\begin{align*}
\ddt \mathcal{I}(\ul{\bm{u}}_\Gamma^{(k)})(0) 
&\approx
\frac{-(1 - c^2) \mathcal{I}(\ul{\bm{u}}_\Gamma^{(k)})(0)
+ \mathcal{I}(\ul{\bm{u}}_\Gamma^{(k)})(\Delta t_0) 
- c^2 \mathcal{I}(\ul{\bm{u}}_\Gamma^{(k)})(\Delta t_0 + \Delta t_1)}
{\Delta t_0 (1 - c)}, 
\\
\dot{\bm{u}}_I^{(1),(k)}(0) 
&\approx
\frac{-(1 - c^2) \bm{u}_I^{(1),(k),0}
+ \bm{u}_I^{(1),(k),1}
- c^2 \bm{u}_I^{(1),(k),2}}
{\Delta t_0 (1 - c)}, \\
c &= \frac{\Delta t_0}{\Delta t_0 + \Delta t_1}.
\end{align*}
Using the SDIRK2 scheme to solve the Neumann problem \eqref{EQ SEMI DISCR NEU}, replacing the heat flux by $\mathcal{I}(\ul{\bm{q}}_1^{(k+1)})$ resp. $\mathcal{I}(\ul{\bm{q}}_2^{(k+1)})$, consists of first solving 
\begin{align*}
\left( \bm{M} + a \Delta t \bm{A} \right) \bm{U}_1^{(2)} 
= \bm{M} \bm{u}^{(k+1), n} - a \Delta t_n
\begin{pmatrix} \bm{0} \\ \mathcal{I}(\ul{\bm{q}}_1^{(k+1)})(t_n + a \Delta t_n)\end{pmatrix},
\end{align*}
followed by
\begin{align*}
& \left( \bm{M} + a \Delta t \bm{A} \right) \bm{u}^{(2), (k+1), n+1} 
= \bm{M} \left(\bm{u}^{(k+1), n} + (1-a) \Delta t_n \bm{k}_1^{(2)}\right) - a \Delta t
\begin{pmatrix} \bm{0} \\ \mathcal{I}(\ul{\bm{q}}_2^{(k+1)})(t_n + \Delta t_n)\end{pmatrix},
\end{align*}
for $n = 0, \ldots, N_2^{(k+1)} - 1$. The relaxation step and termination check are identical to the implicit Euler method. Here, we use a total of only $3$ interpolants: The solution, the heat flux and the stage heat flux. This was achieved by the approximation of $\bm{U}_1^{(2)}\big|_\Gamma$ by $\mathcal{I}(\ul{\bm{u}}_\Gamma^{(k)})$. As can be seen in the numerical results in Section \ref{SEC MR CONV}, this does not lead to a loss of order. Other options for constructing partitioned time-integration methods based on SDIRK2 are discussed in \cite[Chap. 5]{Meisrimel2021}. 

The differences to Algorithm \ref{ALG IE FULL} are in \texttt{SolveDirichlet} and \texttt{SolveNeumann}, which require solving two linear systems each. Additionally, \texttt{SolveDirichlet} now returns two heat fluxes, which are both interpolated and passed into \texttt{SolveNeumann}. 
%
\subsection{Optimal relaxation parameter}\label{SEC THETA}
%
By writing out the linear system for all timesteps in a single WR iteration, one can see that the iteration matrix is a block lower-triangular Toeplitz matrix. Thus, its spectral radius is independent of the number of timesteps, c.f., \cite{Janssen1996}. This means it is sufficient to look at a single timestep to determine the optimal relaxation parameter. In this setting the iteration matrix w.r.t. $\bm{u}_\Gamma^{(k)}$ was already determined in \cite{Monge:2017}. It is given by
\begin{equation*}
\boldsymbol{\Sigma} = -{\bm{S}^{(2)}}^{-1} \bm{S}^{(1)},
\end{equation*}
where
\begin{equation}\label{EQ THETA S MAT}
\bm{S}^{(m)} := \left( \frac{\bm{M}_{\Gamma \Gamma}^{(m)}}{\Delta t} + \bm{A}_{\Gamma \Gamma}^{(m)} \right) - \left( \frac{\bm{M}_{\Gamma I}^{(m)}}{\Delta t} + \bm{A}_{\Gamma I}^{(m)} \right) \left( \frac{\bm{M}_{II}^{(m)}}{\Delta t} + \bm{A}_{II}^{(m)} \right)^{-1} \left( \frac{\bm{M}_{I \Gamma}^{(m)}}{\Delta t} + \bm{A}_{I \Gamma}^{(m)} \right). 
\end{equation}
This is obtained by solving \eqref{EQ IE DIR} for $\bm{u}_I^{(1), (k+1), n+1}$, assuming $\bm{u}_I^{(1),(k+1),n} = \bm{0}$. Insert this result into \eqref{EQ IE FLUX} and then \eqref{EQ IE NEU}. Lastly, solve \eqref{EQ IE NEU} for $\bm{u}_\Gamma^{(k+1),n+1}$, assuming $\bm{u}_{I}^{(2),(k+1),n} = \bm{0}$, using the Schur-complement.

Including the relaxation step yields the following iteration (for a single timestep):
\begin{equation*}
\bm{u}_{\Gamma}^{(k + 1)}
= (\Theta \boldsymbol{\Sigma} + (1 - \Theta) \bm{I}) \bm{u}_{\Gamma}^{(k)}.
\end{equation*}
In the 1D case $\bm{S}^{(m)}$ and $\boldsymbol{\Sigma}$ are scalars, thus the optimal relaxation parameter is
\begin{equation}\label{EQ THETA OPT}
\Theta_{opt} = \frac{1}{\left| 1 + {\bm{S}^{(2)}}^{-1} \bm{S}^{(1)}\right|}.
\end{equation}

In the following we specifically consider a 1D model problem on $\Omega = [-1, 1]$, split at $x_\Gamma = 0$, and an equidistant discretization using linear finite elements. The matrices $\bm{M}^{(m)}_{II}$ and $\bm{A}^{(m)}_{II}$ have a known Toeplitz structure. Thus, one can write down an exact expression of \eqref{EQ THETA S MAT} using an Eigendecomposition to calculate the inverse of $\bm{M}_{II}^{(m)} /\Delta t + \bm{A}_{II}^{(m)}$. Through lengthy but straight forward calculations (see \cite{Monge:2017,Monge2018a}), one obtains the following expressions:
\begin{align}\label{EQ THETA SM}
\begin{split}
 \bm{S}^{(m)}
& = 
\frac{6 \Delta t \Delta x (\alpha_m \Delta x^2 + 3 \lambda_m \Delta t) - (\alpha_m \Delta x^2 - 6 \lambda \Delta t)^2s_m}{18 \Delta t \Delta x^3},\\
s_m
& = 
\sum_{i=1}^{N} \frac{3 \Delta t \Delta x^2 \sin^2 (i \pi \Delta x)}{2 \alpha_m \Delta x^2 + 6 \lambda_m \Delta t + (\alpha_m \Delta x^2 - 6 \lambda_m \Delta t) \cos (i \pi \Delta x)}. 
\end{split}
\end{align}
Using $c = \Delta t/\Delta x^2$, $\Theta_{opt}$ has the following temporal and spatial limits \cite{Monge:2017}: 
\begin{equation}\label{EQ THETA LIMS}
\lim_{c \rightarrow 0} \Theta_{opt} = \frac{\alpha_2}{\alpha_1 + \alpha_2},
\quad 
\lim_{c \rightarrow \infty} \Theta_{opt} = \frac{\lambda_2}{\lambda_1 + \lambda_2}.
\end{equation}
These are consistent with the one-dimensional continuous analysis performed in \cite{gander:16,Mandal:14}. There, a convergence analysis using Laplace transforms for the DNWR method \eqref{EQ CONT DIR PROB}-\eqref{EQ CONT UPDATE PROB} on two identical subdomains $\Omega_1$ and $\Omega_2$ with constant coefficients shows that $\Theta_{opt} = 1/2$. Their result is recovered when approaching the continuous case in the limit $\Delta t/\Delta x^2 \rightarrow \infty$ for constant coefficients.

\begin{figure}
\centering
\includegraphics[width=5.1cm]{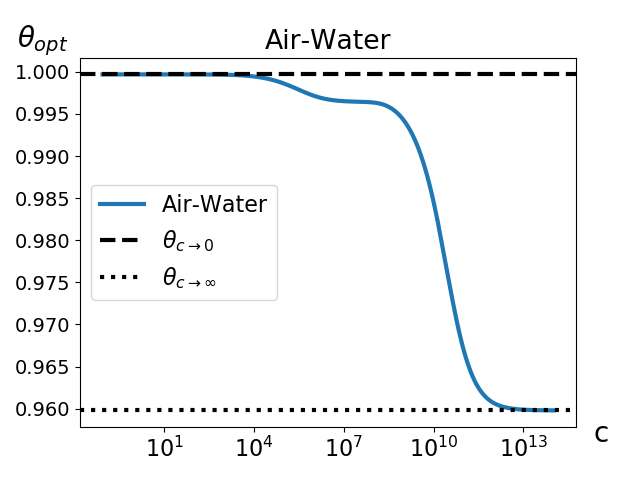} \hfill
\includegraphics[width=5.1cm]{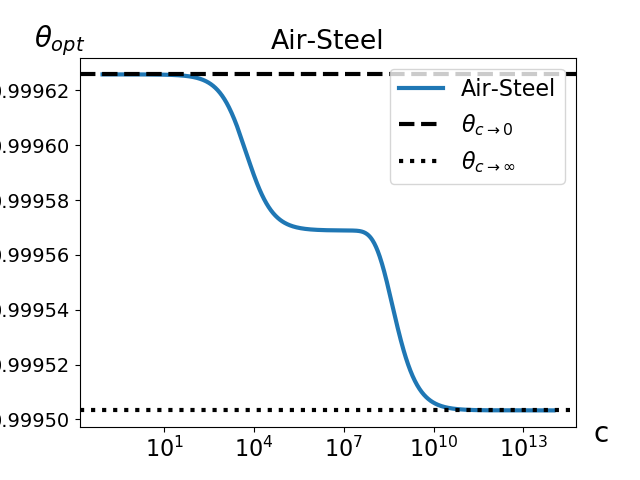} \hfill
\includegraphics[width=5.1cm]{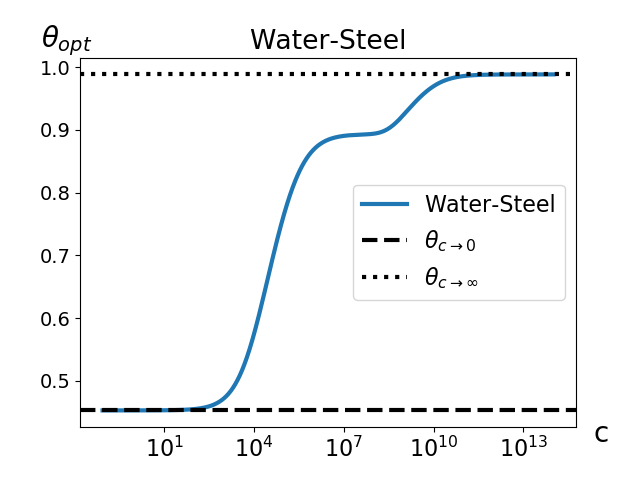} \hfill
\caption{$\Theta_{opt}$ over $\Delta t / \Delta x^2$ for DNWR algorithm. $\Delta x = 1/100$, $T_f = 10^{10}$ and $\Delta t = T_f/2^0, \ldots, T_f/2^{50}$.}
\label{FIG THETA OPT VS CFL}
\end{figure}

Figure \ref{FIG THETA OPT VS CFL} shows $\Theta_{opt}$ for a few material combinations, see Table \ref{TABLE MATERIALS}. One can observe that $\Theta_{opt}$ is continuous and bounded by its spatial and temporal limits \eqref{EQ THETA LIMS}.

%
%
\subsubsection{Multirate relaxation parameter}
%
For $\Delta t_1 \neq \Delta t_2$ the analysis in the previous section does not apply anymore. Instead, we determine $\Theta_{opt}$ based on numerical experiments in Section \ref{SEC THETA MR TEST}. These show that, on average, the optimal choice is to use $\Theta_{opt}$ based on the maximum of $\Delta t_1$ and $\Delta t_2$. This result coincides with the experiments to determine $\Theta_{opt}$ for the multirate NNWR method in \cite{monbirDD25:19}.
%
\section{Time adaptive method}\label{SEC DT CONTROL}
%
The goal of adaptivity is to use step-sizes as large as possible and as small as necessary to reduce computational costs, while ensuring a target accuracy. In particular, using adaptive time stepping for both sub-domains separately, one can attain comparable time-integration errors automatically, bypassing the need to determine a suitable step-size ratio for the multirate case.

The basic idea is to control timestepsizes to keep an error estimate at a given tolerance $TOL$. We use a local error estimate obtained by an embedded technique \cite[chap. IV.8]{HairerII}. With the SDIRK2 method, we obtain a lower order embedded solution $\hat{\bm{u}}_{n+1}$ via
\begin{equation*}
\hat{\bm{u}}_{n+1} = \bm{u}_n + (1 - \hat{a}) \Delta t_n \bm{k}_1 + \hat{a} \Delta t_n \bm{k}_2,
\quad
\hat{a} = 2 - \frac{5}{4} \sqrt{2}.
\end{equation*}
The local error estimate is $\boldsymbol{\ell}_n = \bm{u}_n - \hat{\bm{u}}_n$. We then control timesteps using the proportional-integral (PI) controller
\begin{equation*}
\Delta t_{n+1}
=
\Delta t_n 
\left( \frac{TOL}{\|\boldsymbol{\ell}_n\|_I}\right)^{1/3}
\left( \frac{TOL}{\|\boldsymbol{\ell}_{n-1}\|_I}\right)^{-1/6},
\end{equation*}
cf. \cite{arevalo2017grid}, PI3333. We use this procedure on each subdomain independently. As the initial stepsize we use
\begin{equation*}
\Delta t_0^{(m)}
= 
\frac{T_f \, {TOL^{(m)}}^{1/2}}{100(1 + \| {\bm{M}_{II}^{(m)}}^{-1}\bm{A}_{II}^{(m)} \bm{u}_I^{(m)}(0)\|_I)}, 
\quad m = 1,\,2,
\end{equation*}
c.f. \cite{monbirDD25:19}. We choose the tolerances $TOL^{(m)} = TOL_{WR}/5$, $m = 1,\, 2$. This choice is motivated by \cite{Soderlind06} and already used in a similar context in \cite{BirkenMonge17, monbirDD25:19}. We use the discrete $\mathcal{L}^2$ norm 
\begin{equation}\label{EQ NORM INNER}
\| \bm{u} \|_I^2 = (\bm{u}^T \bm{M} \bm{u})/|\Omega_u|,
\end{equation}
where $\bm{M}$ is the corresponding mass matrix and $|\Omega_u|$ the area on which $\bm{u}$ is defined.

Using this adaptive method, we get independent time-grids for both sub-domains, that are suitable for the given material parameters. A pseudocode of the adaptive SDIRK2 DNWR method is shown in Algorithm \ref{ALG ADAPTIVE}. 
\begin{algorithm}[ht!]
\caption{Pseudocode of the adaptive SDIRK2 DNWR method. The functions \texttt{AdaptiveSolveDirichlet} and \texttt{AdaptiveSolveNeumann} perform the time integration from Section \ref{SEC SDIRK2}, using the step-size control from Section \ref{SEC DT CONTROL}. Note that the time-grids associated with $\ul{\bm{u}}_\Gamma^{(k)}$, $\ul{\bm{q}}_1^{(k+1)}$ and $\ul{\bm{q}}_2^{(k+1)}$ can change in every iteration. This is particularly relevant in the relaxation step, where one needs to interpolate the previous solution to the new time-grid.}
\label{ALG ADAPTIVE}
\begin{algorithmic}
\STATE{\textbf{DNWR\_SDIRK2\_TA}($T_f$, $(\bm{u}_0^{(1)}, \bm{u}_0^{(2)}, \bm{u}_\Gamma(0))$, $\Theta$, $TOL_{WR}$, $k_{\max}$):}
\STATE{$\ul{\bm{u}}_\Gamma^{(0)} \equiv \ul{\bm{u}}_\Gamma(0)$ Initial guess}
\FOR{$k = 0, \ldots, k_{\max} - 1$}
     \STATE{$\mathcal{I}(\ul{\bm{u}}^{(k)}_\Gamma) \gets \text{Interpolation}(\ul{\bm{u}}_\Gamma^{(k)})$}
     \STATE{$\ul{\bm{q}}_1^{(k+1)}, \ul{\bm{q}}_2^{(k+1)}$ $\gets$ \texttt{AdaptiveSolveDirichlet}($T_f$, $TOL/5$, $\bm{u}_0^{(1)}$, $\mathcal{I}(\ul{\bm{u}}^{(k)}_\Gamma)$)}
     \STATE{$\mathcal{I}(\ul{\bm{q}_1}^{(k+1)}) \gets \text{Interpolation}(\ul{\bm{q}}_1^{(k+1)})$}
     \STATE{$\mathcal{I}(\ul{\bm{q}_2}^{(k+1)}) \gets \text{Interpolation}(\ul{\bm{q}}_2^{(k+1)})$}
     \STATE{$\ul{\bm{u}}_\Gamma^{(k+1)}$ $\gets$ \texttt{AdaptiveSolveNeumann}($T_f$, $TOL/5$, $(\bm{u}_0^{(2)}, \bm{u}_0(x_\Gamma))$, $\mathcal{I}(\ul{\bm{q}_1}^{(k+1)})$, $\mathcal{I}(\ul{\bm{q}_2}^{(k+1)})$)}
     \STATE{Compute $\Theta^{(k)}$, see Section \ref{SEC THETA ADAPTIVE}}
    \STATE{$\ul{\bm{u}}_{\Gamma}^{(k+1)} \gets \Theta^{(k)} \ul{\bm{u}}_{\Gamma}^{(k+1)} + (1 - \Theta^{(k)}) \ul{\bm{u}}_{\Gamma}^{(k)}$}
     \IF{$\| \bm{u}_{\Gamma}^{(k+1)}(T_f) - \bm{u}_{\Gamma}^{(k)}(T_f)\|_{\Gamma} < TOL_{WR} \, \| \bm{u}_\Gamma(0)\|_{\Gamma}$}
         \STATE{\textbf{break}}
     \ENDIF
 \ENDFOR
\end{algorithmic}
\end{algorithm}
%
\subsection{Relaxation parameter}\label{SEC THETA ADAPTIVE}
%
In the adaptive method the Dirichlet-Neumann operator changes every WR step, since the timestepsizes change. Consequently, we recompute $\Theta$ in every iteration. We use $\Theta$ as in the multirate case with the average stepsizes from each subdomain. This approach improves upon \cite{monbirDD25:19}, where $\Theta$ was based on the timegrids of the previous WR iteration.
%
\section{The Neumann-Neumann Waveform Relaxation algorithm}
%
Here, we briefly recap the related Neumann-Neumann Waveform relaxation (NNWR) method, c.f. \cite{Fwok:14,MongeBirken:multirate,monbirDD25:19}. Similar to DNWR, NNWR solves \eqref{EQ PROB MONO} in a partitioned manner. The continuous formulation of the algorithm is as follows: Given $g^{(k)}(\bm{x}, t)$ one first solves the following Dirichlet problem on each subdomain:
\begin{align*}
\begin{split}
\alpha_m \frac{\partial u_m^{(k+1)}(\bm{x},t)}{\partial t} - \nabla \cdot (\lambda_m \nabla u_m^{(k+1)}(\bm{x},t)) &= 0, \quad (\bm{x}, t) \in \Omega_m \times (0, T_f],\\
u_m^{(k+1)}(\bm{x},t) &= 0, \quad (\bm{x},t) \in \partial \Omega_m \backslash \Gamma \times [0, T_f], \\
u_m^{(k+1)}(\bm{x},t) &= g^{(k)}(\bm{x},t), \quad (\bm{x},t) \in \Gamma \times (0, T_f], \\
u_m^{(k+1)}(\bm{x},0) &= u_m^{0}(\bm{x}), \quad \bm{x} \in \Omega_m.
\end{split}
\end{align*}
Next, one solves the following Neumann problems:
\begin{align}
\begin{split}\label{EQ NNWR NEU PROB}
\alpha_m \frac{\partial \psi_m^{(k+1)}(\bm{x},t)}{\partial t} - \nabla \cdot (\lambda_m \nabla \psi_m^{(k+1)}(\bm{x},t)) &= 0, \quad (\bm{x},t) \in \Omega_m \times (0, T_f],\\
\psi_m^{(k+1)}(\bm{x},t) &= 0, \quad (\bm{x},t) \in \partial \Omega_m \backslash \Gamma \times [0, T_f], \\
\lambda_m \frac{\partial \psi_m^{(k+1)}(\bm{x},t)}{\partial \bm{n}_m} 
&= \lambda_1 \frac{\partial u_1^{(k+1)}(\bm{x},t)}{\partial \bm{n}_1}
+ \lambda_2 \frac{\partial u_2^{(k+1)}(\bm{x},t)}{\partial \bm{n}_2}, \quad (\bm{x},t) \in \Gamma \times (0, T_f], \\
\psi_m^{(k+1)}(\bm{x},0) &= 0, \quad \bm{x} \in \Omega_m.
\end{split}
\end{align}
Finally, the update step is
\begin{equation*}
g^{(k+1)}(\bm{x}, t) = g^{(k)}(\bm{x}, t) - \Theta(\psi^{(k+1)}_1(\bm{x}, t) + \psi^{(k+1)}_2(\bm{x}, t)), \quad (\bm{x},t) \in \Gamma \times [0, T_f].
\end{equation*}
For the fully discrete version and a detailed algorithmic description see \cite{MongeBirken:multirate}.

One can solve on the subdomains in parallel. The NNWR algorithm is based on the exact same Dirichlet and Neumann problems as the DNWR algorithm, but with different input data, namely fluxes and initial value for the Neumann problem \eqref{EQ NNWR NEU PROB}. Hence one can directly use the time-discretizations as described in Section \ref{SEC TIME DISCR}, including time-adaptivity, c.f. \cite{monbirDD25:19}. 

Under the same restrictions as in Section \ref{SEC THETA}, i.e., $\Omega = [-1,1]$, split at $x_\Gamma = 0$ and linear finite elements on an equidistant grid, one can analogously calculate $\Theta_{opt}$:
\begin{equation}\label{EQ NNWR THETA OPT}
\Theta_{opt} = \frac{1}{|2 + {\bm{S}^{(1)}}^{-1}\bm{S}^{(2)} + {\bm{S}^{(2)}}^{-1}\bm{S}^{(1)}|}
\end{equation}
with $\bm{S}^{(m)}$ given by \eqref{EQ THETA SM}, see \cite{MongeBirken:multirate} for details. The spatial and temporal limits based on $c = \Delta t/\Delta x^2$ are
\begin{equation}\label{EQ NNWR THETA LIMS}
\lim_{c \rightarrow 0} \Theta_{opt} = \frac{\alpha_1 \alpha_2}{(\alpha_1 + \alpha_2)^2},
\quad 
\lim_{c \rightarrow \infty} \Theta_{opt} = \frac{\lambda_1 \lambda_2}{(\lambda_1 + \lambda_2)^2}.
\end{equation}
%
%
\section{Numerical results}
%
We now present numerical experiments to illustrate the validity of the theoretical results and to test the robustness of the relaxation parameters in 2D, SDIRK2 and with multirate resp. adaptive time-grids. The methods and algorithms described have been implemented in Python 3.6, the code is available at \cite{DNWR_code}.

We consider the domains $\Omega = [-1, 1]$ for 1D and $\Omega = [-1, 1] \times [0, 1]$ for 2D, with $\Omega_1$ and $\Omega_2$ split at $x_\Gamma = 0$. Our initial conditions are
\begin{align}\label{EQ U0 1}
u(x) = 500 \sin((x+1) \pi/2),
\quad \text{resp.} \quad
u(x, y) = 500 \sin((x+1) \pi/2) \sin(y \pi).
\end{align}
As the coefficients $\alpha$ and $\lambda$ in \eqref{EQ PROB MONO} we consider the materials as shown in Table \ref{TABLE MATERIALS}.
\begin{table}[ht!]
\begin{center}
\begin{tabular}{|c|c|c|}
\hline \textbf{Material} & $\alpha = \rho \cdot c_p [J/(K m^3)]$ & $\lambda [W/(m K)]$ \\
\hline Air & $1.293 \cdot 1005$
 & $0.0243$
 \\
\hline Water & $999.7 \cdot 4192.1
$ & $0.58

$ \\
\hline Steel & $7836 \cdot 443$
 & $48.9$ \\
\hline 
\end{tabular}
\caption{Material parameters.
\label{TABLE MATERIALS}}
\end{center}
\end{table}

The resulting heterogeneous cases are Air-Water, Air-Steel and Water-Steel. We use $T_f = 10^4$ in all cases.

As space discretization we use linear finite elements on equidistant grids (1D: $\Delta x = 1/200$, 2D: $\Delta x = 1/100$) as shown in Figure \ref{FIG DOMAIN FE}, see \cite{BirkenMonge17} for more details. The resulting linear equation systems are solved with direct solvers.

We define our multirate setup via $N$ as the number of base timesteps. In subdomain $\Omega_m$ we then use $N_m = c_m \cdot N$ timesteps. We consider the following cases: Coarse-coarse ($c_1 = c_2 = 1$), coarse-fine ($c_1 = 1$, $c_2 = 10$) and fine-coarse ($c_1 = 10$, $c_2 = 1$).
%
\subsection{Multirate relaxation parameter}\label{SEC THETA MR TEST}
%
The question is what $\Theta$ to choose for DNWR in the multirate case, i.e., which $\Delta t$ to use in \eqref{EQ THETA SM}. We consider the following four choices: \textbf{Max}/\textbf{Min}/\textbf{Avg} by taking the maximum, minimum or average of $\Delta t_1$ and $\Delta t_2$ to compute $\Theta_{opt}$, and "\textbf{Mix}": $\Theta_{opt} = 1/(| 1 + {\bm{S}^{(2)}}^{-1}(\Delta t_2) \bm{S}^{(1)}(\Delta t_1)|)$.

We experimentally determine the convergence rate via
\begin{equation}
\| \bm{u}_\Gamma^{(k)}(T_f) - \bm{u}_\Gamma^{(k-1)}(T_f)\|_{\Gamma}\, /\, \|\bm{u}_\Gamma^{(k-1)}(T_f) - \bm{u}_\Gamma^{(k-2)}(T_f)\|_{\Gamma},
\end{equation}
i.e., the reduction rate in the update. Here, we perform up to $k_{\max} = 6$ iterations and take the mean of the update reductions, but never the last iteration, which could be near machine precision. This experiment is done using IE for the 1D test case and $N=1$, as we aim to determine the asymptotic convergence rates.
\begin{figure}[ht!]
\centering
\includegraphics[width=7cm]{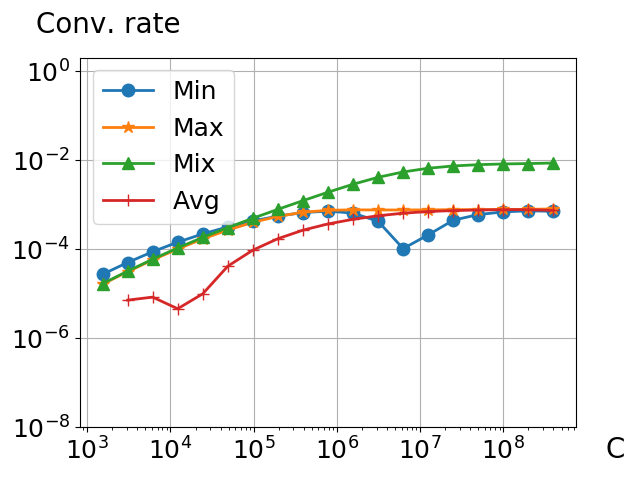} \hfill
\includegraphics[width=7cm]{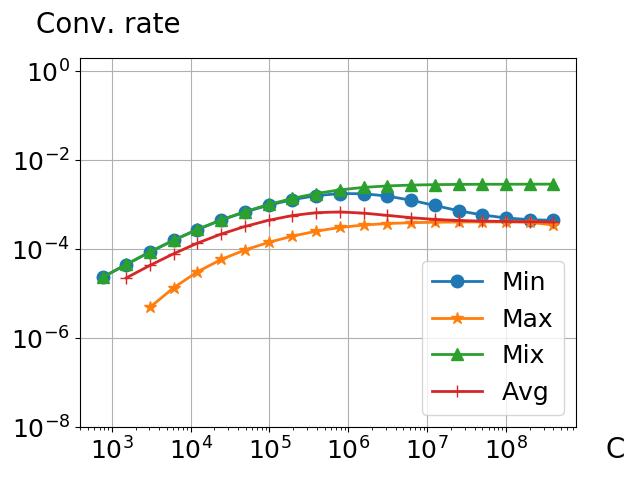} \hfill
\caption{Observed convergence rates over $c = \Delta t / \Delta x^2$ for DNWR IE, 1D, air-water. Left: Coarse-fine. Right: Fine-coarse.}
\label{FIG THETA MR AIR WATER}
\end{figure}
\begin{figure}[ht!]
\centering
\includegraphics[width=7cm]{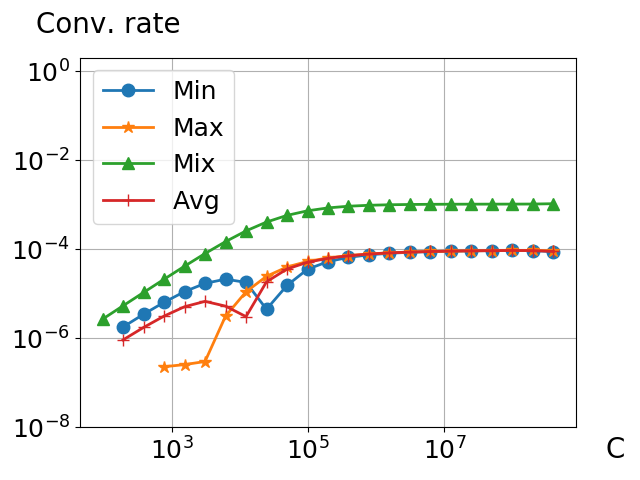} \hfill
\includegraphics[width=7cm]{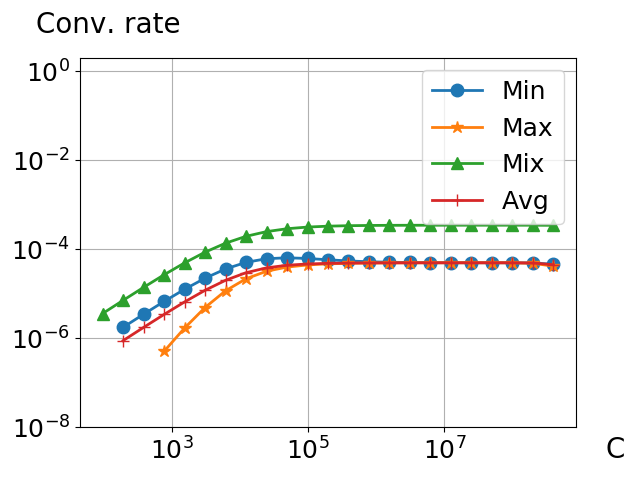} \hfill
\caption{Observed convergence rates over $c = \Delta t / \Delta x^2$ for DNWR IE, 1D, air-steel. Left: Coarse-fine. Right: Fine-coarse.}
\label{FIG THETA MR AIR STEEL}
\end{figure}
\begin{figure}[ht!]
\centering
\includegraphics[width=7cm]{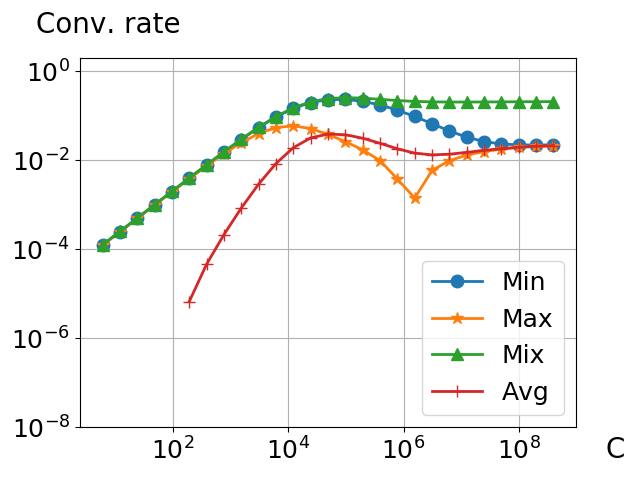} \hfill
\includegraphics[width=7cm]{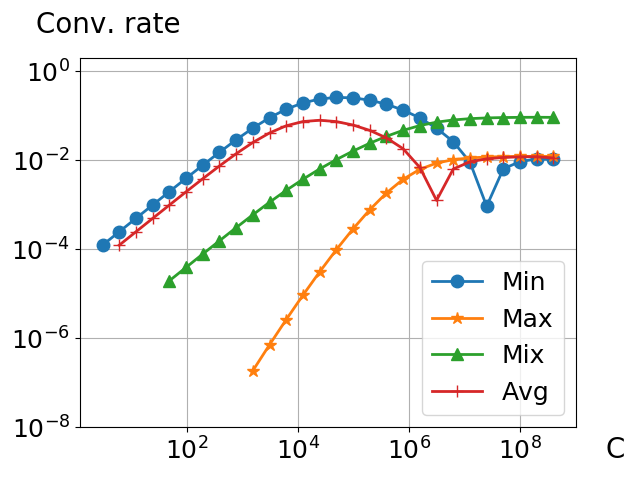} \hfill
\caption{Observed convergence rates over $c = \Delta t / \Delta x^2$ for DNWR IE, 1D, water-steel. Left: Coarse-fine. Right: Fine-coarse.}
\label{FIG THETA MR WATER STEEL}
\end{figure}

The results in Figures \ref{FIG THETA MR AIR WATER}, \ref{FIG THETA MR AIR STEEL} and \ref{FIG THETA MR WATER STEEL} show that all options yield comparable results, with "\textbf{Max}" being most consistent, making it our choice for DNWR in the multirate setting.

Numerical experiments in \cite{MongeBirken:multirate} yielded the same conclusion for the NNWR method. As such we use \eqref{EQ NNWR THETA OPT} based on the larger step-size.
%
\subsection{Optimality of relaxation parameter}
%
We now verify the optimality of the relaxation parameters \eqref{EQ THETA OPT} for DNWR and \eqref{EQ NNWR THETA OPT} for NNWR in 1D with implicit Euler and also test the convergence rate and robustness. To this end, we determine the experimental convergences rates as in Section \ref{SEC THETA MR TEST}, for varying $\Theta$ in 1D and 2D for both implicit Euler and SDIRK2. Time-integration is done in multirate and non-multirate, up to $T_f$ using $N = 100$ base timesteps.

We expect little variation for implicit Euler and SDIRK2, since SDIRK2 consists of two successive implicit Euler steps. We anticipate more notable differences between 1D and 2D. Lastly, convergence rates might deviate due to transitive effects of WR, since the iteration matrices are non-normal.

In the plots, the blue highlighted range on the $x$-axis marks the spatial and temporal limits of $\Theta$, see \eqref{EQ THETA LIMS} resp. \eqref{EQ NNWR THETA LIMS}.
%
\subsubsection{DNWR}
%
\begin{figure}[ht!]
\centering
\includegraphics[width=5.1cm]{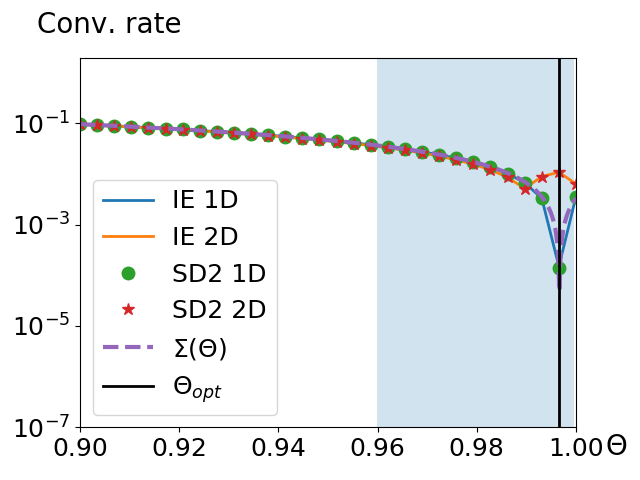} \hfill
\includegraphics[width=5.1cm]{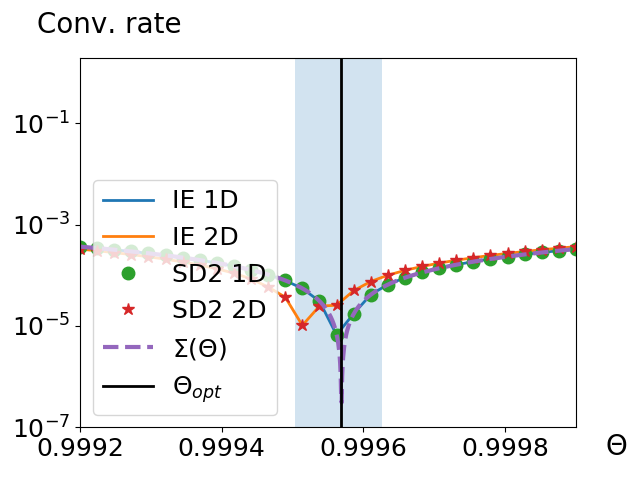} \hfill
\includegraphics[width=5.1cm]{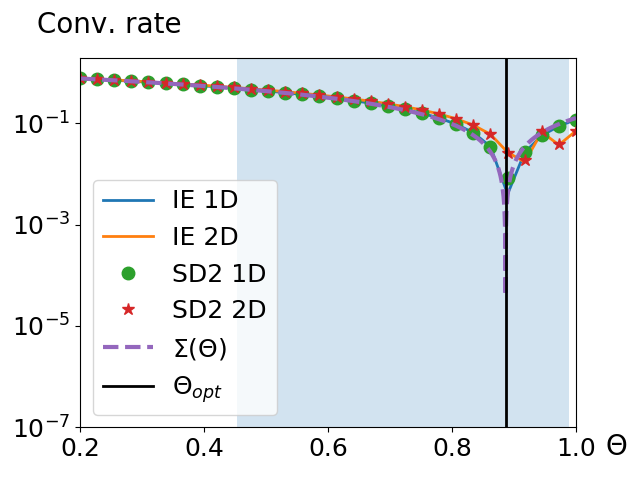} \hfill
\caption{Left: Air-water, fine-coarse. Centre: Air-steel, coarse-coarse. Right: Water-steel, coarse-fine. Observed convergence rates for DNWR algorithm.}
\label{FIG OBS CONV}
\end{figure}
Results are seen in Figure \ref{FIG OBS CONV}. In both 1D and 2D, SDIRK2 rates match those of implicit Euler. In all cases, 1D results closely align with the theoretical result marked by $\Sigma(\Theta)$. 2D results are slightly off, but still yield good error reduction rates using the 1D $\Theta_{opt}$. For air-water, the error reduction rate is $\approx 10^{-2}$, i.e., the coupling residual gains two decimals in accuracy per iteration. The air-steel coupling yields very fast convergence with an error reduction rate of $\approx 10^{-4}$ and water-steel rates are between $0.1$ and $0.01$ for $\Theta_{opt}$.

Additionally, we see that DNWR is convergent for all shown $\Theta$. Thus, in the worst case DNWR convergence is slow, yet not divergent.
%
\subsubsection{DNWR - non-square geometry}
\begin{figure}[ht!]
\centering
\includegraphics[width=5.1cm]{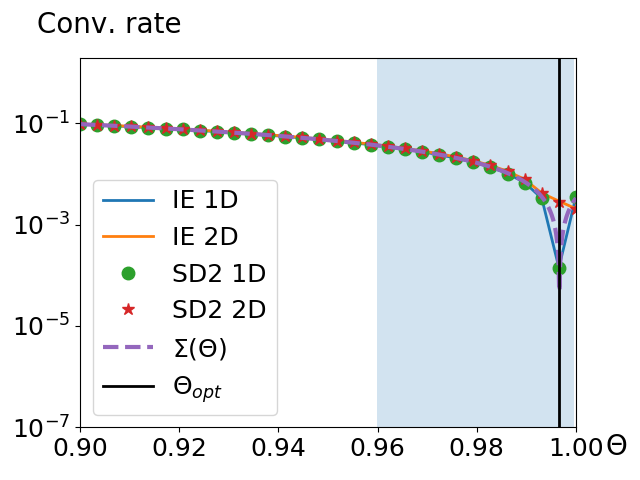} \hfill
\includegraphics[width=5.1cm]{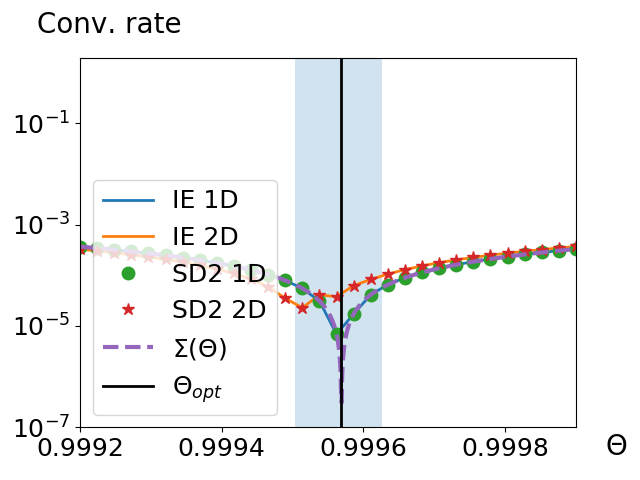} \hfill
\includegraphics[width=5.1cm]{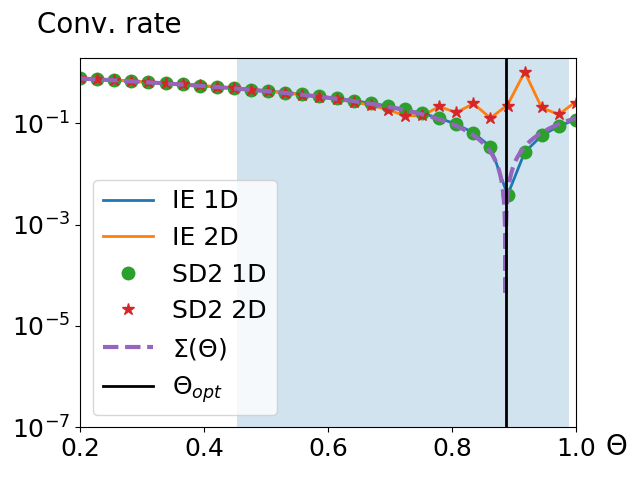} \hfill
\caption{Left: Air-water. Centre: Air-steel. Right: Water-steel. Observed convergence rates for DNWR algorithm using non-square geometry and matching step-sizes.}
\label{FIG OBS CONV NON SQUARE}
\end{figure}

%
We test if the DNWR results extend to non-square domains. In particular, we consider the spatial domain with $x \in [-9, 1]$, $x_{\Gamma} = 0$. We use the initial conditions 
\begin{align*}
u(x) = 500 \sin((x+9) \pi/10),
\quad \text{resp.} \quad
u(x, y) = 500 \sin((x+9) \pi/10) \sin(y \pi).
\end{align*}
We test only the non-multirate setting. Results are shown in Figure \ref{FIG OBS CONV NON SQUARE} and strongly resemble the results for square, identical domains in Figure \ref{FIG OBS CONV}, except for slower convergence rates for the 2D water-steel case.

This similarity can be explained by looking at the Schur-complements \eqref{EQ THETA S MAT}. $\bm{M}_{\Gamma \Gamma}^{(m)}$ and $\bm{A}_{\Gamma \Gamma}^{(m)}$ remain unchanged. $\bm{M}_{II}^{(1)}$ resp. $\bm{A}_{II}^{(1)}$ increase in size, but values and structure persist. The matrices $\bm{M}_{I\Gamma}^{(1)}$, $\bm{A}_{I\Gamma}^{(1)}$ and $\bm{M}_{\Gamma I}^{(1)}$, $\bm{A}_{\Gamma I}^{(1)}$ are padded with additional zeros. The increased size of $\bm{M}_{II}^{(1)}/\Delta t + \bm{A}_{II}^{(1)}$ does affect its inverse, but after multiplication with $\bm{M}_{I\Gamma}^{(1)}/\Delta t + \bm{A}_{I\Gamma}^{(1)}$ and $\bm{M}_{\Gamma I}^{(1)}/\Delta t + \bm{A}_{\Gamma I}^{(1)}$, which are mostly zero, the effect on $\bm{S}^{(1)}$ is expected to be minor.
%
\subsubsection{NNWR}
%
\begin{figure}[ht!]
\centering
\includegraphics[width=5.1cm]{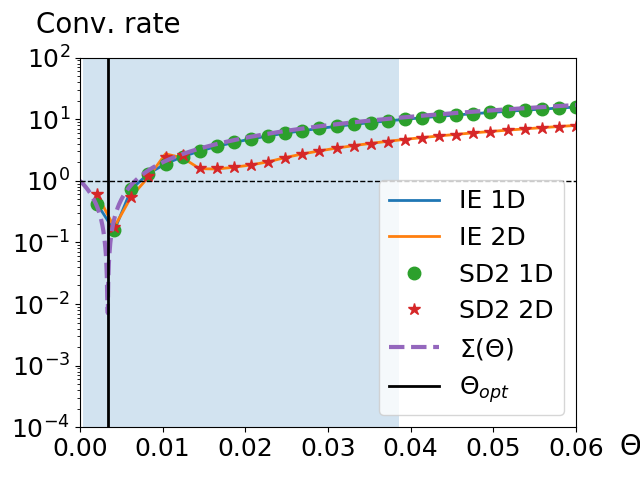}\hfill
\includegraphics[width=5.1cm]{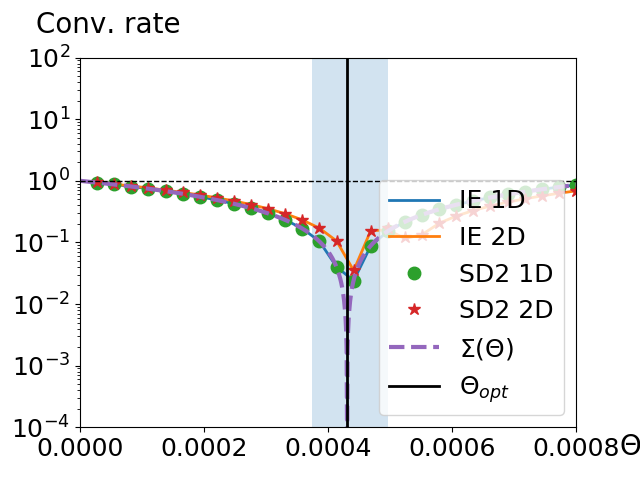}\hfill
\includegraphics[width=5.1cm]{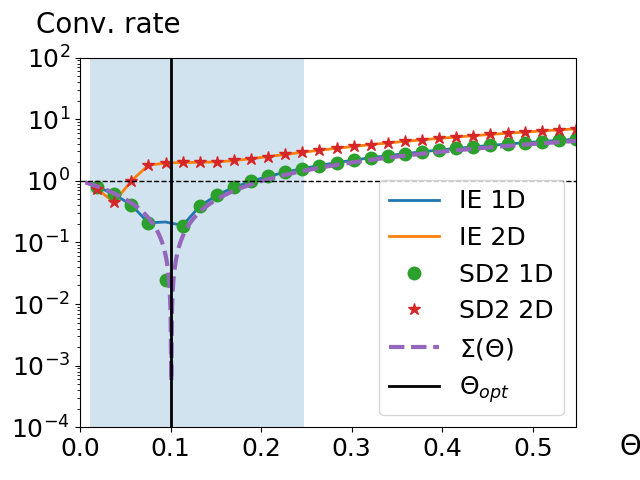}\hfill
\caption{Left: Air-water, fine-coarse. Centre: Air-steel, coarse-coarse. Right: Water-steel, coarse-fine. Observed convergence rates for NNWR algorithm.}
\label{FIG NNWR OBS CONV}
\end{figure}
Results are shown in Figure \ref{FIG NNWR OBS CONV}. We additionally mark the divergence limit at $1$, showing the range of viable $\Theta$ is very small for NNWR and unlike DNWR, relaxation is non-optional for convergence. In particular, one may get divergence for $\Theta$ within the range marked by the temporal and spatial limits. Convergence rates for implicit Euler and SDIRK2 results are almost identical. {1D convergence rates align well with the theoretical results in all cases. In 2D, the air-water and air-steel results match with the 1D results, yielding rates of about $0.01 - 0.1$. However, water-steel shows divergence in 2D, when using $\Theta_{opt}$. In the convergent cases, the observed error reduction rates are slower than for DNWR, {this is particularly pronounced in the air-steel case, with a difference of about 3 orders of magnitude. Overall, NNWR shows a lack of robustness.

One might achieve better convergence rates using macrostepping, i.e., successively performing the algorithm on smaller time-windows. This may speed up convergence on each time-window, but the coupling residual propagates through erroneous initial values. On the other hand, DNWR performs well on the given time-windows.
%
\subsection{Multirate - convergence order of time-integration}\label{SEC MR CONV}
%
We show convergence of the error, on the whole domain in the discrete $\mathcal{L}^2$ norm \eqref{EQ NORM INNER} and using $T_f = 1$, for $\Delta t \rightarrow 0$. Our reference solution is the monolithic solution for sufficiently small step-sizes, thus measuring both the time-integration error and coupling residual.

Results for $TOL_{WR} = 10^{-13}$ can be seen in Figure \ref{FIG DNWR MR CONV} for DNWR and in Figure~\ref{FIG NNWR MR CONV} for NNWR. We attain the expected first and second order convergence rates for $\Delta t \rightarrow 0$.
\begin{figure}[h!]
\centering
\includegraphics[width=5.1cm]{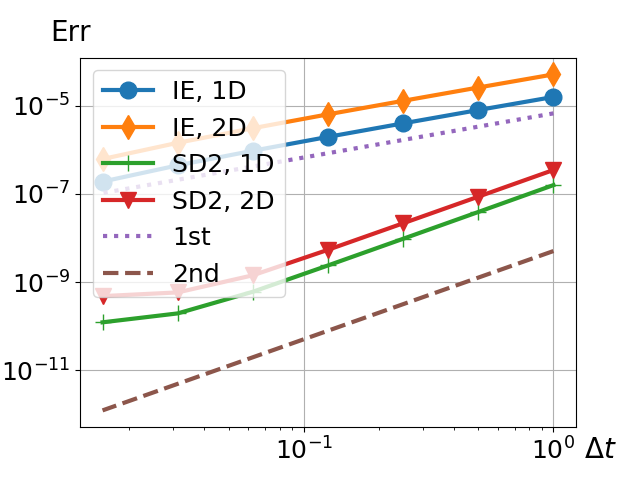} \hfill
\includegraphics[width=5.1cm]{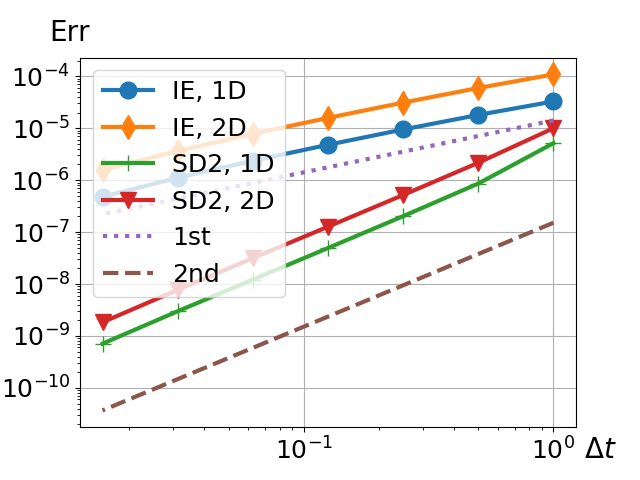} \hfill
\includegraphics[width=5.1cm]{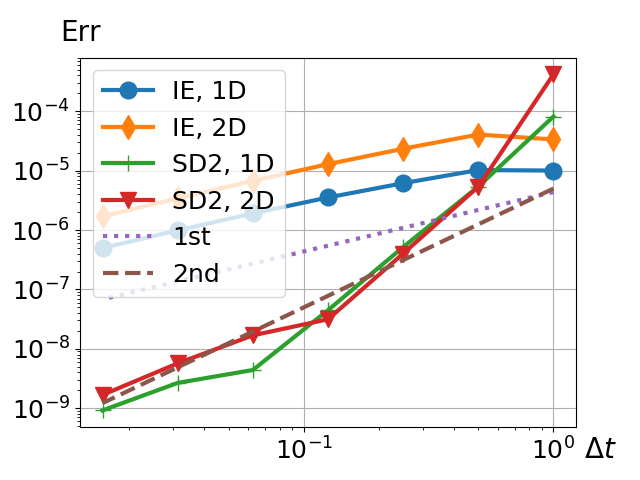} \hfill
\caption{Left: Air-water, fine-coarse. Centre: Air-steel, coarse-coarse. Right: Water-steel, coarse-fine. Error over $\Delta t$ for DNWR and $T_f = 1$.}
\label{FIG DNWR MR CONV}
\end{figure}

\begin{figure}[h!]
\centering
\includegraphics[width=5.1cm]{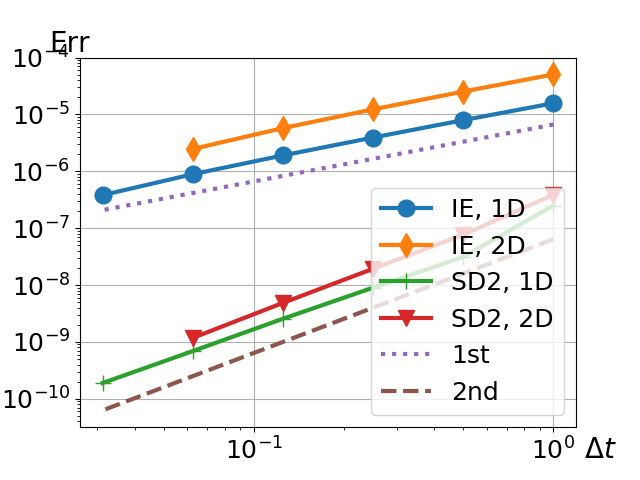} \hfill
\includegraphics[width=5.1cm]{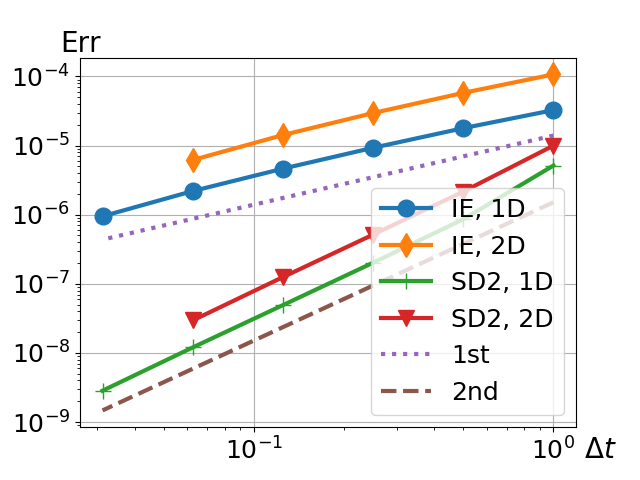} \hfill
\includegraphics[width=5.1cm]{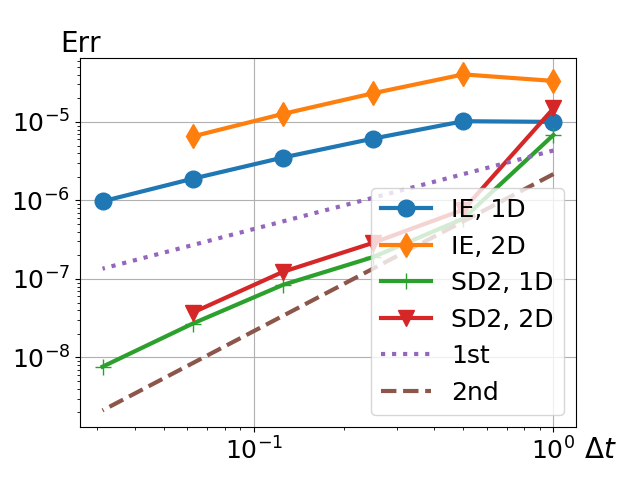} \hfill
\caption{Left: Air-water, fine-coarse. Centre: Air-steel, coarse-coarse. Right: Water-steel, coarse-fine. Error over $\Delta t$ for NNWR and $T_f = 1$.}
\label{FIG NNWR MR CONV}
\end{figure}
%
\subsection{Time adaptive results}
%
We consider the time-adaptive DNWR method described in Section \ref{SEC DT CONTROL}. The reference for error computation is the solution using $TOL_{WR} = 10^{-8}$ in 1D and $TOL_{WR} = 10^{-7}$ in 2D. We expect the errors to be proportional to the tolerance for $TOL_{WR} \rightarrow 0$, which is observed in Figure \ref{FIG TA}. Due to its lacking robustness, we do not consider time-adaptive NNWR.
\begin{figure}[h!]
\centering
\includegraphics[width=7cm]{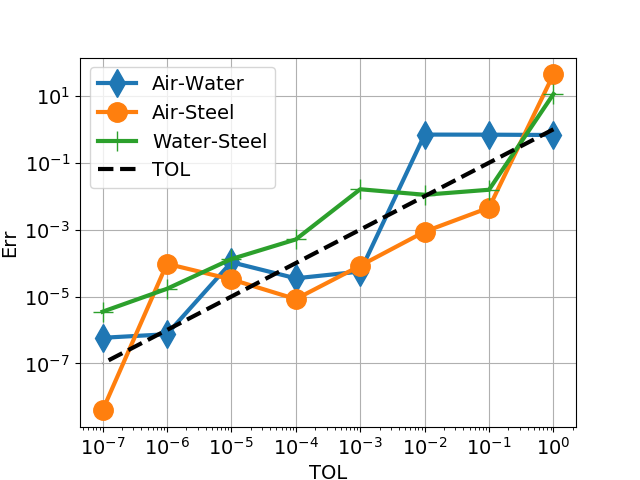} \hfill
\includegraphics[width=7cm]{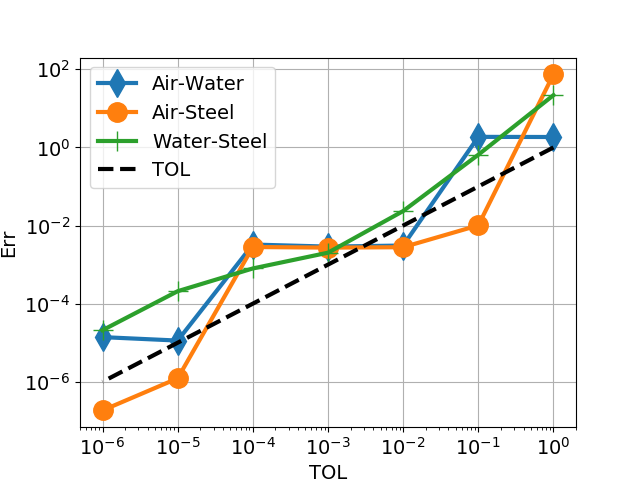} \hfill
\caption{Left: 1D. Right: 2D. Error over $TOL_{WR}$ for the time-adaptive DNWR method.}
\label{FIG TA}
\end{figure}
%
\subsubsection{Error over work comparison}
%
We now compare efficiency of the adaptive and multirate method for the 2D test case with $\Delta x = 1/200$. For this we compare error over work, which we measure as the total number of timesteps.

We choose the stepsize ratios in the multirate setting such that both domains use comparable \textit{CFL} numbers, which is achieved by $c_2 = c_1 D_2/D_1$, $D_m = \lambda_m/\alpha_m$, $m = 1, 2$. However, we require $c_m \in \mathbb{N}$. W.l.o.g., assume $D_2/D_1 > 1$, we then set $c_1 = 1$ and round down $c_2 = D_2/D_1$. See Table \ref{TABLE ERR WORK STEPS} for the resulting stepsize ratios for our material configurations.

To compare the multirate method with the time-adaptive method, we parametrize the former by the number of base timesteps $N$. Given $\Delta t_m = T_f/(c_m \cdot N)$, $m = 1,\,2$, we compute the associated time integration error $e_{\Delta t_1, \Delta t_2}$ using $TOL_{WR} = 10^{-12}$ and a monolithic reference solution with $\Delta t = \min(\Delta t_1, \Delta t_2)/2$. We then use $TOL_{WR} = e_{\Delta t_1, \Delta t_2}/5$ in the termination criterion for the multirate method, for its error over work comparison. Finally, our references for the error computations in the error over work comparison are adaptive solutions with $TOL_{WR} = 10^{-6}$. 

Results are shown in Figure \ref{FIG WORK ERROR 1} with the resulting stepsizes ratios for the adaptive case in Table \ref{TABLE ERR WORK STEPS}. The adaptive method is 4 times more efficient in the water-steel case, of similar efficiency in the air-water case and less efficient in the air-steel case. This can be explained by the stepsize ratios in Table \ref{TABLE ERR WORK STEPS}. The closer the multirate stepsize ratios correspond to the adaptive ones, the better the performance of multirate in comparison with the adaptive method.
\begin{table}[ht!]
\begin{center}
\begin{tabular}{|c||c|c|c|}
\hline
& \textbf{Air-water} & \textbf{Air-steel} & \textbf{Water-steel} \\
\hline
multirate ($c_1:c_2$) & $135:1$ & $1:1$ & $1:101$\\
\hline
adaptive $u_0^1$ & $33.88:1$ & $1.15:1$& $1:2.47$\\
\hline
adaptive $u_0^2$ & $21.25:1$& $1.09:1$ & $1:3.20$\\
\hline
\end{tabular}
\caption{Timestep ratios for the multirate and adaptive method (final grid) by materials. $u_0^1$ is for the initial condition \eqref{EQ U0 1} and $u_0^2$ is for \eqref{EQ U0 2}. 
\label{TABLE ERR WORK STEPS}}
\end{center}
\end{table}
\begin{figure}[ht!]
\centering
\includegraphics[width=5.1cm]{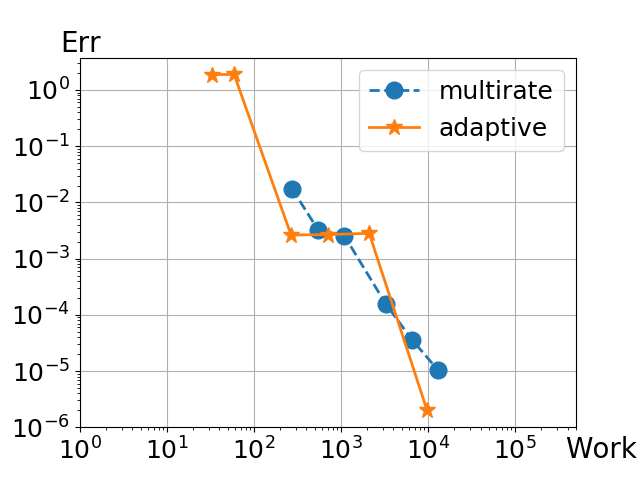} \hfill
\includegraphics[width=5.1cm]{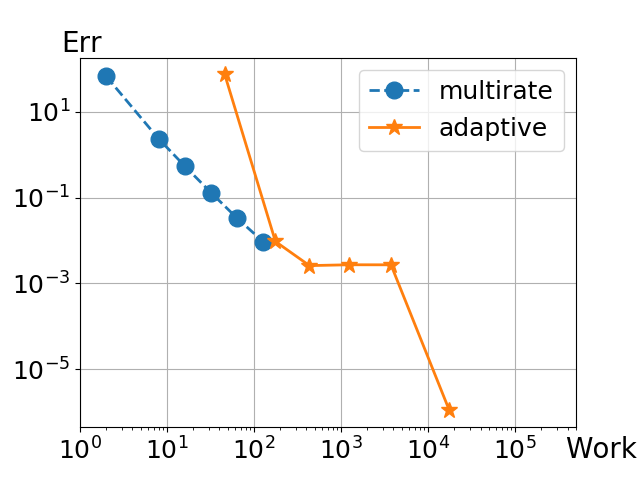} \hfill
\includegraphics[width=5.1cm]{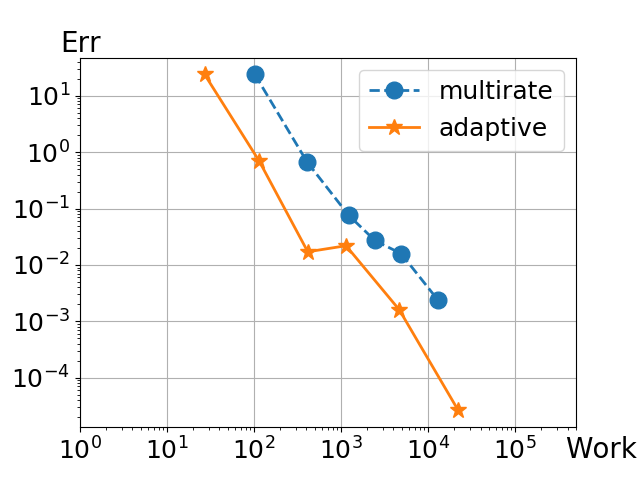} \hfill
\caption{Left: Air-water. Centre: Air-steel. Right: Water-steel. DNWR work over error comparison for 2D test case using initial condition \eqref{EQ U0 1}.}
\label{FIG WORK ERROR 1}
\end{figure}

As second test case we consider the initial condition
\begin{align}\label{EQ U0 2}
u(x, y) = 800 \sin((x+1) \pi)^2 \sin(y \pi).
\end{align}
Here, we have $\bm{u}_{\Gamma}(0) = \bm{0}$ and thus skip the relative norm for the termination check \eqref{EQ TERMINATION CRIT}. Results are shown in Figure \ref{FIG WORK ERROR 2} with stepsize ratios in Table \ref{TABLE ERR WORK STEPS}. In the air-steel case performance is approximately equal, whereas adaptive performance is about 4 resp. 25 times better in the air-water resp. water-steel case.

Overall we see that performance depends on the stepsize ratios. This makes the adaptive method a more robust choice, since it automatically determines suitable stepsize ratios, which vary for e.g., different initial conditions.
\begin{figure}[ht!]
\centering
\includegraphics[width=5.1cm]{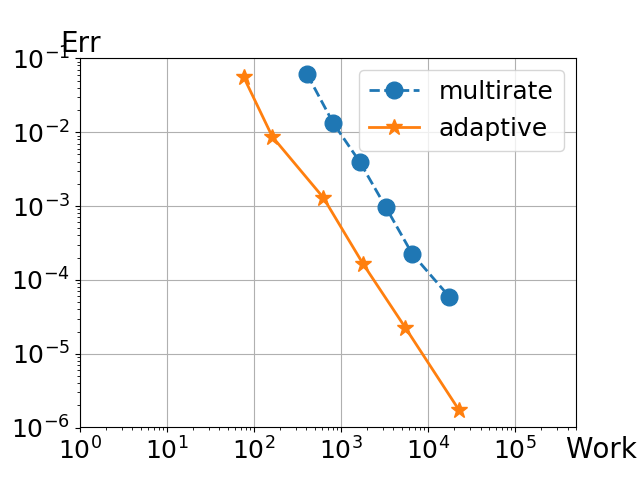} \hfill
\includegraphics[width=5.1cm]{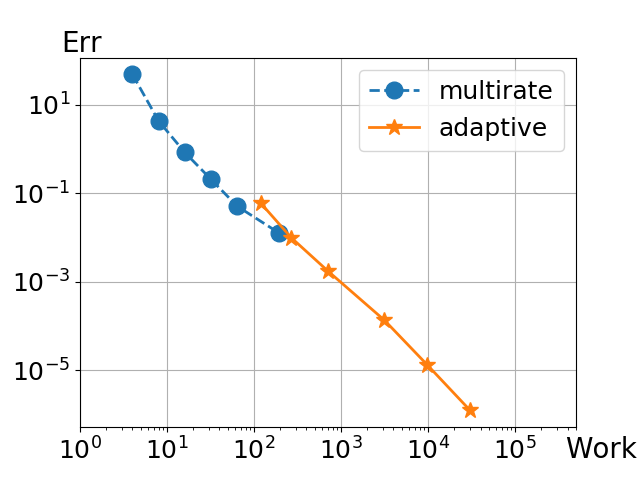} \hfill
\includegraphics[width=5.1cm]{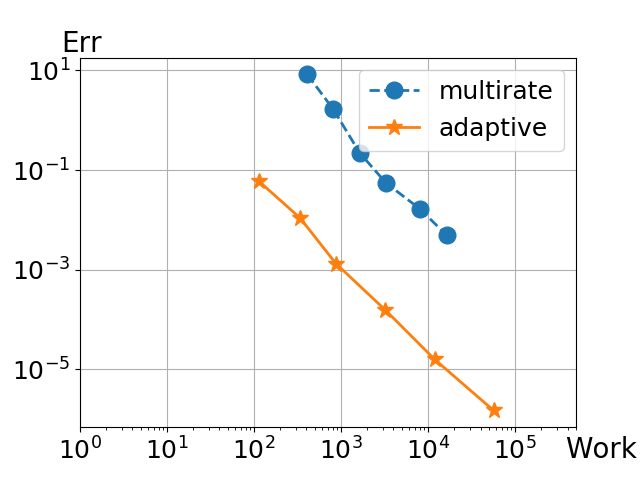} \hfill
\caption{Left: Air-water. Centre: Air-steel. Right: Water-steel. DNWR work over error comparison for 2D test case using initial condition \eqref{EQ U0 2}.}
\label{FIG WORK ERROR 2}
\end{figure}
%
\section{Summary and conclusions} 
%
We derived first and second order, multirate resp. time-adaptive DNWR methods for heterogeneous coupled heat equations. The optimal relaxation parameter $\Theta_{opt}$ for WR is shown to be identical to the one for the basic DN iteration. We experimentally show how to adapt $\Theta$ in the multirate case. The observed convergence rates using an analytical $\Theta_{opt}$ for 1D implicit Euler are shown to be very robust, yielding fast convergence rates for a second order method and 2D, for various material combinations and multirate settings on long time intervals.

The same tests for the related NNWR methods employing identical Dirichlet and Neumann subsolvers, using an analytical $\Theta_{opt}$ for 1D implicit Euler, show a lack of robustness, possibly resulting in divergence.

The time-adaptive DNWR method is experimentally shown to be favorable over mutirate, due ease of use and superior performance. The latter is due to the resulting stepsizes being more suitably chosen than those of the multirate solver. Overall, we obtain a fast, robust, time adaptive (on each domain), partitioned solver for unsteady conjugate heat transfer.
\bibliographystyle{siam}
\bibliography{literature}
\end{document}